\documentclass[10pt]{article}
\usepackage{amsmath}
\usepackage{amssymb}
\makeatletter
\@addtoreset{equation}{section}
\makeatother

 \newtheorem{thm}{Theorem}[section]
 \newtheorem{cor}[thm]{Corollary}
 \newtheorem{lem}[thm]{Lemma}
 \newtheorem{prop}[thm]{Proposition}
\newtheorem{defn}[thm]{Definition}
\newtheorem{rem}[thm]{Remark}
 \numberwithin{equation}{section}
\newtheorem{hyp}[thm]{Hypothesis}
\def\bR {\mathbf{R}}

\def\l {{\lambda}}

\begin{document}

\title{Bursting Dynamics of the 3D Euler Equations  in Cylindrical Domains}

%----------Author 1
\author{Fran\c cois Golse~\footnote{golse@math.polytechnique.fr}
\,\,\,\footnote{and Laboratoire J.-L. Lions, Universit\'e Paris Diderot-Paris 7} \\
   Ecole Polytechnique, CMLS\\
   91128 Palaiseau Cedex, France\\
   \mbox{   } \\
 Alex Mahalov~\footnote{mahalov@asu.edu}
  and Basil Nicolaenko~\footnote{byn@stokes.la.asu.edu}\\
   Department of Mathematics and Statistics\\
   Arizona State University\\
   Tempe, AZ 85287-1804, USA}

%%% ----------------------------------------------------------------------
\date{  }
\maketitle

\begin{abstract}
A class of three-dimensional initial data 
characterized by uniformly large vorticity
is considered for the 3D incompressible Euler equations
in bounded cylindrical domains.
The  fast singular oscillating limits
of the 3D Euler equations are investigated for parametrically 
resonant cylinders. Resonances of fast oscillating swirling
Beltrami waves deplete the Euler nonlinearity. 
These waves  are exact solutions of the 3D Euler equations.
We construct the 3D resonant Euler systems; the latter are 
countable uncoupled and coupled $SO(3;\mathbf{C})$ and $SO(3; \mathbf{R})$
rigid body systems. They conserve both energy and helicity.
The 3D resonant Euler systems are vested with bursting dynamics,
where the ratio of the enstrophy at time $t=t^*$ to the enstrophy
at $t=0$ of some remarkable orbits becomes very large for very small 
times $t^*$; similarly for higher norms $\mathbf{H}^s,\ s\geq 2$.
These orbits are topologically close to homoclinic cycles. For the time 
intervals where  $\mathbf{H}^s$ norms, $s\geq 7/2$ of the limit resonant
orbits do not blow up, we prove that the full 3D Euler equations possess 
smooth solutions
close to the resonant orbits  uniformly in strong norms.
\end{abstract}



%-------------------------------------------------
  \section{Introduction}
%-------------------------------------------------

The issues of blowup of smooth solutions  and finite time singularities of
the vorticity field for 3D incompressible Euler equations are still a major 
open problem. The Cauchy problem in 3D bounded axisymmetric cylindrical 
domains is attracting considerable attention:
with bounded, smooth, \textit{non}-axisymmetric 3D initial data, under
the constraints of conservation of \textit{bounded} energy, can the 
vorticity field blow up in \textit{finite} time?
Outstanding numerical claims for this have recently been disproven \cite{Ke},
\cite{Hou1}, \cite{Hou2}.
The classical analytical criterion  of Beale-Kato-Majda \cite{B-K-M} for
non-blow up in finite time requires the time integrability of the $L^{\infty}$ norm of the vorticity. 
DiPerna and Lions \cite{Li} have given examples of global weak solutions of the 3D Euler equations
which are smooth (hence unique) if the initial conditions are smooth (specifically in
$\mathbf{W}^{1,p}(D),\ p>1$).
However, these flows are really 2-Dimensional in $x_1, x_2$, 3-components flows,
independent from the third coordinate $x_3$. Their examples  \cite{DiPe-Li}
show that   solutions (even smooth ones) of the 3D Euler equations cannot be estimated
in $\mathbf{W}^{1,p}$ for $1 < p<\infty$ on any time interval $(0,\, T)$ 
if the initial data are only assumed to be bounded in $\mathbf{W}^{1,p}$.
   Classical local existence theorems in 3D 
bounded or periodic domains by Kato \cite{Ka}, 
Bourguignon-Br\'{e}zis \cite{Bou-Br}
and Yudovich \cite{Yu1}, \cite{Yu2} require
some minimal smoothness for the initial conditions (IC), e.g., in $\mathbf{H}^s(D),\ s>\frac{5}{2}$.

The classical formulation for the Euler equations is
\begin{eqnarray}
\partial_t \mbox{\bf{V}} + (\mbox{\bf{V}} \cdot \nabla ) \mbox{\bf{V}}
 =
- \nabla p, \mbox{   } \nabla \cdot \mbox{\bf{V}}=0, \hspace{.2in} & & \label{nsch8}\\
 \mbox{\bf{V}}\cdot \mathbf{N} = 0 \mbox{ on } \partial D ,\hspace{.85in}
 & & \label{nsch8/bc}
\end{eqnarray}
 where $\partial D$ is the boundary of a bounded, connected domain $D$, $\mathbf{N}$ 
the normal to $\partial D$, $\mathbf{V}(t,y)=(V_1, V_2, V_3)$  the velocity field,
$y=(y_1, y_2, y_3)$, and $p$ is the pressure.

The equivalent Lam\'{e} form \cite{Ar-Khe}
\begin{eqnarray}
  \partial_t\mathbf{V}+ curl\mathbf{V}\times\mathbf{V} + 
    \nabla \left( p+\frac{1}{2}|\mathbf{V}|^2\right) = 0, && \label{Lame_form}\\
   \nabla\cdot \mathbf{V} = 0, \hspace{1in}&& \label{Lame_form/div}
\end{eqnarray}
\begin{subequations}\label{Lame_vorticity_new}
   \begin{align}
  \partial_t \boldsymbol{\omega}+ curl(\boldsymbol{\omega}\times \mathbf{V} ) = 0, 
   \label{Lame_vorticity} \\
   \boldsymbol{\omega}= curl\mathbf{V}, \ \ \ \ 
\end{align}
\end{subequations}
implies  conservation of Energy:
\begin{equation}
  E(t) = \frac{1}{2}\int_D |\mathbf{V}(t,y)|^2\,dy.
\end{equation}
The helicity $Hel(t)$\, \cite{Ar-Khe}, \cite{Mof}, is  conserved:
\begin{equation}
  Hel(t) = \int_D \mathbf{V}\cdot\boldsymbol{\omega}\,dy,
\end{equation}
for $D = \mathbf{R}^3$ and when $D$ is a periodic lattice.
Helicity is also conserved for cylindrical domains,
provided that $\boldsymbol{\omega}\cdot\mathbf{N}=0$ on the cylinder's
lateral boundary at $t=0$ (see \cite{M-N-B-G}).

From the theoretical point of view,
the principal difficulty in the analysis of 3D Euler equations
is due to the presence of the vortex stretching term
$(\boldsymbol{\omega}\cdot \nabla )\mathbf{V}$ in the vorticity equation 
(\ref{Lame_vorticity}).
The equations (\ref{Lame_form}) and (\ref{Lame_vorticity}) are equivalent to:
\begin{equation}\label{Lie}
  \partial_t\boldsymbol{\omega} + [\boldsymbol{\omega}, \mathbf{V}] = 0, 
\end{equation}
where $[a, b] = \mbox{ curl }(\mathbf{a}\times\mathbf{b})$ is the commutator
in the infinite dimensional Lie algebra of divergence-free vector fields \cite{Ar-Khe}.
This point of view has led to celebrated developments in Topological Methods
in Hydrodynamics \cite{Ar-Khe}, \cite{Mof}.
The striking analogy between the Euler equations for hydrodynamics and the
Euler equations for a rigid body (the latter associated to the Lie Algebra of the Lie group
$SO(3,\mathbf{R}))$ had already been pointed out by Moreau \cite{Mor1};
Moreau was the first to demonstrate conservation of Helicity (1961) \cite{Mor2}.
This has led to extensive speculations to what extent/in what cases are the 
solutions of the 3D Euler equations ``close'' to those of coupled 3D rigid body equations in
some asymptotic sense. Recall that the Euler equations for a rigid body in $\mathbf{R}^3$ 
is:
\begin{subequations}\label{rigid}
  \begin{align}
  \mathbf{m}_t+\boldsymbol{\omega}\times\mathbf{m} = 0, \ \
  \mathbf{m} = A\boldsymbol{\omega}, \label{rigid_a}\\
  \mathbf{m}_t+[\boldsymbol{\omega}, \mathbf{m}] = 0, \hspace{.3in} \label{rigid_b} 
  \end{align}
\end{subequations}
where $\mathbf{m}$ is the vector of angular momentum relative to the body,
$\boldsymbol{\omega}$ the angular velocity in the body and $A$ the inertia operator 
\cite{Ar1}, \cite{Ar-Khe}.

The Russian school of Gledzer, Dolzhansky, Obukhov \cite{G-D-O} and Vishik 
\cite{Vish}
has extensively investigated dynamical systems of hydrodynamic type and their applications. 
They have considered hydrodynamical models built upon generalized rigid body systems 
in $SO(n,\mathbf{R})$, following Manakhov \cite{Man}. Inspired by turbulence physics,
they have investigated ``shell'' dynamical systems modeling turbulence cascades; 
albeit such systems are flawed as they only preserve energy, not helicity.
To address this, they have constructed and studied in depth $n$-dimensional dynamical systems with quadratic homogeneous nonlinearities and {\bf two} quadratic first integrals 
$F_1, F_2$. Such systems can be written using sums of Poisson brackets:
\begin{equation}\label{eq1.10}
  \frac{dx^{i_1}}{dt} = \frac{1}{2}\sum_{i_2,...,i_n}\epsilon^{i_1i_2...i_n}
    p_{i_4...i_n}\left( \frac{\partial F_1}{\partial x^{i_2}} 
                        \frac{\partial F_2}{\partial x^{i_3}} -
                        \frac{\partial F_1}{\partial x^{i_3}}
                        \frac{\partial F_2}{\partial x^{i_2}}\right) ,
\end{equation}
where constants $p_{i_4...i_n}$ are antisymmetric in $i_4,..., i_n$.

A simple version of such a quadratic hydrodynamic system was introduced by Gledzer \cite{Gl1} in 1973. A deep open issue of the work by the Gledzer-Obukhov school is
whether there exist indeed classes of I.C. for the 3D Cauchy Euler problem 
(\ref{nsch8}) for which solutions are actually asymptotically close in strong norm, on arbitrary large time intervals to solutions of such hydrodynamic systems, with conservation of 
both energy and helicity. Another unresolved issue is the blowup or global regularity 
for the ``enstrophy'' of such systems when their dimension $n\rightarrow \infty$.

This article reviews some current new results of a research program in the spirit of the
Gledzer-Obukhov school; this program builds-up on the results of
\cite{M-N-B-G} for 3D Euler in bounded cylindrical domains.
Following the original approach of \cite{B-M-N1}-\cite{B-M-N4} in periodic domains,
\cite{M-N-B-G} prove the 
non blowup of the 3D incompressible
Euler equations for a class of
three-dimensional initial data characterized by uniformly large vorticity
in bounded cylindrical domains.
There are no conditional assumptions on the
properties of solutions at later times, nor are the global
solutions close to some 2D manifold.
The initial vortex stretching is large.
The approach of
proving regularity is based on investigation of fast
singular oscillating limits and nonlinear averaging methods
in the context of almost periodic functions \cite{Bo-Mi}, \cite{Bes},
\cite{Cor}.
Harmonic analysis tools based on curl eigenfunctions and eigenvalues are crucial.
One establishes the global regularity of the 3D limit resonant Euler
equations without any restriction on the size of 3D initial data.
The resonant Euler equations are characterized by a depleted nonlinearity. 
After establishing strong convergence to the limit resonant
equations, one
bootstraps this into the regularity on arbitrary large time
intervals of the solutions of 3D Euler Equations with weakly
aligned uniformly large vorticity at $t=0$. 
\cite{M-N-B-G} theorems hold for generic cylindrical domains, 
for a set of height/radius ratios of full Lebesgue measure. For such cylinders, 
the 3D limit resonant Euler equations are restricted to two-wave resonances of 
the vorticity waves and are vested with an infinite countable number of
new conservation laws.
The latter are adiabatic invariants for the original 3D Euler equations.

Three-wave resonances exist for a nonempty countable set of $h/R$ 
($h$ height, $R$ radius of the cylinder) and   moreover accumulate in the limit
of vanishingly small vertical (axial) scales.
This is akin to Arnold tongues~\cite{Ar2} for the Mathieu-Hill equations  and
raises nontrivial issues of possible singularities/lack thereof for dynamics ruled by 
infinitely many resonant triads at vanishingly small axial scales.
In such a context, the 3D resonant Euler equations do conserve the energy and helicity of
the field.

In this review, we consider cylindrical domains with
parametric resonances in $h/R$ and investigate
in depth the structure and dynamics of 3D resonant Euler systems. These
 parametric 
resonances in $h/R$ are proven to be non-empty.
Solutions to Euler equations with uniformly large initial vorticity are
expanded along a full \textit{complete} basis of elementary
swirling waves ($\mathbf{T}^2$ in time). Each such quasiperiodic, dispersive
vorticity wave is a quasiperiodic Beltrami flow; these are exact solutions
of 3D Euler equations with vorticity parallel
to velocity. There are no Galerkin-like truncations in the decomposition of
the full 3D Euler field. The Euler equations, restricted to resonant triplets of these
dispersive Beltrami waves, determine the ``resonant Euler systems''.
The basic ``building block'' of these (a priori $\infty$-dimensional)
systems are proven to be $SO(3;\mathbf{C})$ and $SO(3;\mathbf{R})$ rigid body systems:
  \begin{equation}
\label{eq1.1}
\begin{aligned}
\dot{U}_k+(\l_m-\l_n)U_mU_n&=0
\\
\dot{U}_m+(\l_n-\l_k)U_nU_k&=0
\\
\dot{U}_n+(\l_k-\l_m)U_kU_m&=0
\end{aligned}
\end{equation}
These $\lambda$'s are \textit{eigenvalues of the curl operator}
in the cylinder,  $curl\mathbf{\Phi}_n^{\pm}=\pm\lambda_n\mathbf{\Phi}_n^{\pm}$;
the curl eigenfunctions are steady elementary Beltrami flows, and the dispersive
Beltrami  waves oscillate
with the frequencies
$\pm \frac{h}{2\pi \epsilon}\frac{n_3}{\lambda_n},\ n_3$ vertical
wave number (vertical shear), $0<\epsilon <1$.
Physicists \cite{Ch-Ch-Ey-H} have computationally
demonstrated the physical impact of the polarization of Beltrami modes 
$\mathbf{\Phi}^\pm$ on intermittency in the joint cascade of energy
and helicity in turbulence. 

Another ``building block'' for resonant Euler systems is a pair of
$SO(3;\mathbf{C})$ or $SO(3;\mathbf{R})$ rigid bodies coupled
via a common principal axis of inertia/mo- ment of inertia:
\begin{subequations}\label{eq1.2}
   \begin{align}
     \dot{a}_k&=(\lambda_m-\lambda_n)  \Gamma a_m a_n \\
      \dot{a}_m&=(\lambda_n-\lambda_k)  \Gamma a_n a_k \\
     \dot{a}_n&=(\lambda_k-\lambda_m)  \Gamma a_k a_m 
            + (\lambda_{\tilde{k}}-\lambda_{\tilde{m}})\tilde{\Gamma}a_{\tilde{k}}a_{\tilde{m}}\\
     \dot{a}_{\tilde{m}}&=(\lambda_n-\lambda_{\tilde{k}})\tilde{\Gamma}a_n a_{\tilde{k}} \\
     \dot{a}_{\tilde{k}}&=(\lambda_{\tilde{m}}-\lambda_n)\tilde{\Gamma}a_{\tilde{m}} a_n,
   \end{align}
\end{subequations}
where $\Gamma$ and $\tilde{\Gamma}$ are parameters in $\mathbf{R}$ defined in
Theorem \ref{Theorem4.10}.
Both resonant systems (\ref{eq1.1}) and (\ref{eq1.2})
 conserve  energy and helicity. We prove that the dynamics of these resonant
systems admit equivariant families
of homoclinic cycles connecting hyperbolic critical points.
We demonstrate bursting dynamics: the ratio
\[ ||\mathbf{u}(t)||_{H^s}^2/||\mathbf{u}(0)||_{H^s}^2,\ s\geq 1
\] 
can burst arbitrarily large on arbitrarily small times, for properly
chosen parametric domain resonances $h/R$.
Here 
\begin{equation}\label{Hs}
||\mathbf{u}(t)||_{H^s}^2 = \sum_n (\lambda_n)^{2s} |u_n(t)|^2\,.
\end{equation}
The case $s=1$ is the enstrophy. The ``bursting'' orbits are topologically
close to the homoclinic cycles.

Are such dynamics for the resonant systems relevant to the full 3D 
Euler equations (\ref{nsch8})-(\ref{Lie})? The answer lies in the following crucial
``shadowing'' Theorem \ref{theorem2.10}.
Given the same initial conditions, given the maximal  time interval
$0\leq t< T_m$ where the resonant orbits of the resonant Euler equations 
do {\bf not} blow up, then the \textit{strong} norm $\mathbf{H}^s$ of the difference
between the exact Euler orbit and the resonant orbit is uniformly small
on $0\leq t< T_m$, 
provided that the vorticity of the I.C. is large enough.
Paradoxically, the larger the vortex streching of the I.C., 
the better the uniform approximation. This deep result is based on cancellation of 
fast oscillations in \textit{strong} norms, in the context of almost periodic
functions of time with values in Banach spaces
(Section 4 of \cite{M-N-B-G}).  It includes uniform approximation in the spaces
$\mathbf{H}^s,\ s>5/2$.
For instance, given a quasiperiodic orbit on some time torus $\mathbf{T}^l$
for the resonant Euler systems, the exact solutions to the Euler equations will
remain $\epsilon$-close to the resonant
quasiperiodic orbit on a time interval $0\leq t\leq \max T_i$,
$1\leq i\leq l$, $T_i$ elementary periods, for large enough initial vorticity.
If orbits of the resonant Euler systems admit bursting dynamics in the strong
norms $\mathbf{H}^s,\ s\geq 7/2$, so do
\textit{some exact solutions of the full 3D Euler equations}, for properly
chosen parametrically resonant cylinders.

%-------------------------------------------------------------------
  \section{Vorticity waves and resonances of elementary swirling flows}
%-------------------------------------------------------------------

We study initial value problem for the three-dimensional Euler
equations with initial data characterized by uniformly large
vorticity:
\begin{eqnarray}
\label{nsch8_a}
\partial_t \mbox{\bf{V}} + (\mbox{\bf{V}} \cdot \nabla) \mbox{\bf{V}}
 =
- \nabla p, \mbox{   } \nabla \cdot \mbox{\bf{V}}=0, \hspace{.2in} && \\
\label{nsch8_a/ic} \mbox{\bf{V}}(t,y)|_{t=0} =
\mbox{\bf{V}}(0)=\tilde{\mbox{\bf{V}}}_0(y) + \frac{\Omega}{2} \mathbf{e}_3
\times y && 
\end{eqnarray}
where $y=(y_1,y_2,y_3)$, $\mbox{\bf{V}}(t,y)=(V_1, V_2, V_3)$ is
the velocity field and $p$ is the pressure. In Eqs.~(\ref{nsch8})
$\mathbf{e}_3$ denotes the vertical unit vector and $\Omega$ is a constant
parameter. The field $\tilde{\mbox{\bf{V}}}_0(y)$ depends on three
variables $y_1$, $y_2$ and $y_3$. Since $curl (\frac{\Omega}{2}
\mathbf{e}_3 \times y) = \Omega \mathbf{e}_3$, the vorticity vector at initial time
$t=0$ is
\begin{equation}
\label{curl0} curl \mbox{\bf{V}}(0,y)=curl
\tilde{\mbox{\bf{V}}}_0(y) + \Omega \mathbf{e}_3,
\end{equation}
and the initial vorticity has a large component weakly aligned
along $\mathbf{e}_3$, when $\Omega >> 1$. These are fully three-dimensional
large initial data with large initial 3D vortex stretching.
We denote by $\mathbf{H}_{\sigma}^s$ the usual Sobolev space of solenoidal vector fields. 

The base flow
\begin{equation}
  \mathbf{V}_s(y) = \frac{\Omega}{2}\mathbf{e}_3\times y, \ curl \mathbf{V}_s(y)
=\Omega \mathbf{e}_3
\end{equation}
is called a steady swirling flow and is a steady state solution
(\ref{nsch8})-(\ref{Lame_form/div}), as  $curl(\Omega \mathbf{e}_3\times \mathbf{V}_s(y)) = 0$.
In  (\ref{nsch8_a/ic}) and (\ref{curl0}), we consider I.C. which are an arbitrary 
{\bf\textit{(not small)}} perturbation of the base swirling flow $\mathbf{V}_s(y)$ and introduce
\begin{eqnarray}
  \mathbf{V}(t,y) = \frac{\Omega}{2}\mathbf{e}_3\times y + \mathbf{\tilde{V}}(t,y), \hspace{1.5in} && \\
  curl \mathbf{V}(t,y) = \Omega \mathbf{e}_3+curl \mathbf{\tilde{V}}(t,y), \hspace{1.3in}  && \\
  \partial_t \mathbf{\tilde{V}}+ curl\mathbf{\tilde{V}}\times\mathbf{\tilde{V}} + 
  \Omega \mathbf{e}_3\times \mathbf{\tilde{V}} +
   curl \mathbf{\tilde{V}}\times \mathbf{V}_s(y) +\nabla p^\prime = 0,\ \nabla\cdot\mathbf{\tilde{V}}=0, \label{swiring} && \\
  \mathbf{\tilde{V}}(t,y)|_{t=0}=\tilde{\mathbf{V}}_0(y).\hspace{2in} \label{eq2.8}&& 
\end{eqnarray}

 Eqs.~(\ref{nsch8_a}) and (\ref{swiring}) are studied in cylindrical domains
\begin{equation}
\label{cylinder} \mbox{\bf{C}}=\{ (y_1, y_2, y_3) \in
\mbox{\bf{R}}^3: 0 < y_3 < 2 \pi/\alpha, \mbox{   } y_1^2 + y_2^2
< R^2 \}
\end{equation}
where $\alpha$ and $R$ are positive real numbers. If $h$ is
the height of the cylinder, $\alpha = 2\pi /h$. Let
\begin{equation}
\label{cylb} \Gamma=\{ (y_1, y_2, y_3) \in \mbox{\bf{R}}^3: 0 <
y_3 < 2 \pi/\alpha, \mbox{   } y_1^2 + y_2^2 = R^2 \}.
\end{equation}
Without loss of generality, we can assume that $R=1$.
Eqs.~(\ref{nsch8_a}) are considered with periodic boundary
conditions in $y_3$
\begin{equation}
\label{perbc}
\mbox{\bf{V}}(y_1,y_2,y_3)=\mbox{\bf{V}}(y_1,y_2,y_3+2 \pi/\alpha)
\end{equation}
and vanishing normal component of velocity on $\Gamma$
\begin{equation}\label{nbc} 
  \mbox{\bf{V}} \cdot {\bf N}=\tilde{{\bf V}} \cdot {\bf
N}=0 \mbox{   on }  \Gamma ;
\end{equation}
where ${\bf N}$ is the normal vector to $\Gamma$. From the
invariance of 3D Euler equations under the symmetry $y_3
\rightarrow -y_3$, $V_1 \rightarrow V_1$, $V_2 \rightarrow V_2$,
$V_3 \rightarrow - V_3$, all results in this article extend to
cylindrical domains bounded by two horizontal plates. Then the
boundary conditions in the vertical direction are zero flux on the
vertical boundaries (zero vertical velocity on the plates). One
only needs to restrict vector fields to be even in $y_3$ for
$V_1$, $V_2$ and odd in $y_3$ for $V_3$, and double the
cylindrical domain to $-h \leq y_3 \leq + h$.

We choose $\tilde{\mbox{\bf{V}}}_0(y)$ in ${\bf H}^s(\mathbf{C}), \ s>5/2$.
In \cite{M-N-B-G},  for the case of ``non-resonant cylinders'', 
that is, non-resonant $\alpha = 2\pi/h$, we have 
established regularity for arbitrarily large finite times
for the 3D Euler solutions for $\Omega$ large, but {\it finite}.
Our solutions are not close in any sense to those of the 2D or
``quasi 2D'' Euler and they are characterized by fast oscillations
in the $\mathbf{e}_3$ direction, together with a large vortex stretching
term 
\[ \boldsymbol{\omega}(t,y) \cdot \nabla {\bf V}(t,y)= \omega_1
\frac{\partial \mathbf{V}}{\partial y_1} + \omega_2 \frac{\partial
\mathbf{V}}{\partial y_2} + \omega_3 \frac{\partial \mathbf{V}}{\partial y_3},\ \ 
t \geq 0
\]
 with leading component $\Big| \Omega \frac{\partial}{\partial
y_3} \mathbf{V}(t,y)\Big| \gg 1$. There are no assumptions on oscillations in
$y_1$, $y_2$ for our solutions (nor for the initial condition
$\tilde{\mbox{\bf{V}}}_0(y))$.

Our approach is entirely based on sturying fast singular oscillating limits
of Eqs.~(\ref{nsch8})-(\ref{Lame_vorticity}), nonlinear averaging and
cancelation of oscillations in the nonlinear interactions for the
vorticity field for large $\Omega$. This has been developed in
\cite{B-M-N2}, \cite{B-M-N3}, and \cite{B-M-N4} for
the cases of periodic lattice domains and the infinite space
$\mathbf{R}^3$. 

It is well known that fully three-dimensional 
initial conditions with uniformly large vorticity excite
fast Poincar\'{e} vorticity waves \cite{B-M-N2}, \cite{B-M-N3}, 
\cite{B-M-N4},
\cite{Poi}.
Since individual
Poincar\'{e} wave modes are related to the eigenfunctions
of the curl operator, they are exact
time-dependent solutions of the full
nonlinear 3D Euler equations. Of course, their linear superposition
does not preserve this property. Expanding solutions of 
(\ref{nsch8_a})-(\ref{eq2.8}) along such vorticity waves demonstrates potential 
nonlinear resonances of such waves. First recall spectral properties of
the $curl$ operator in bounded, connected domains:

\begin{prop} {\em (\cite{M-N-B-G}) } The curl operator admits a self-adjoint 
extension under the zero flux boundary conditions, with a discrete real
spectrum $\lambda_n=\pm |\lambda_n|,\ |\lambda_n|>0$ for every $n$
and $|\lambda_n|\rightarrow +\infty $ as $|n|\rightarrow \infty$.
The corresponding eigenfunctions $\boldsymbol{\Phi}_n^{\pm}$
\begin{equation}
   curl \boldsymbol{\Phi}_n^{\pm} = \pm |\lambda_n| \boldsymbol{\Phi}_n^{\pm}
\end{equation}
are complete in the space
\begin{equation}
   \boldsymbol{J}^0 = \left\{ \boldsymbol{U} \in L^2(D):\,
     \nabla\cdot\mathbf{U} =0  \mbox{ and } \mathbf{U}\cdot\mathbf{N}|_{\partial D}=0
    \mbox{ and } \int_o^h\mathbf{U}\,dz = 0 \right\}.
\end{equation}
\end{prop}

\begin{rem} In cylindrical domains, with cylindrical coordinates 
$(r, \theta , z)$, the eigenfunctions admit the representation:
\begin{equation}
   \boldsymbol{\Phi}_{n_1,n_2,n_3} = 
       \left( \Phi_{r,n_1,n_2,n_3}(r), \Phi_{\theta ,n_1,n_2,n_3}(r), 
                 \Phi_{z,n_1,n_2,n_3}(r) \right) e^{in_2\theta}e^{i\alpha n_3z},
\end{equation}
with $n_2= 0, \pm 1, \pm 2,...$, $n_3=\pm 1, \pm 2, ...$ and $n_1=0, 1, 2,...$. Here $n_1$ indexes the
eigenvalues of the equivalent Sturm-Liouville problem in the radial coordinates,
and $n=(n_1,n_2,n_3)$.  
See {\em \cite{M-N-B-G}} for technical details. From now on,
we use the generic variable $z$ for any vertical (axial) coordinate $y_3$ or $x_3$. 
  For $n_3=0$ (vertical averaging along the axis of the cylinder), 
  2-Dimensional, 3-component solenoidal fields must be 
  expanded along a complete basis for fields derived from 2D stream functions: 
 \begin{eqnarray*}
   \mathbf{\Phi}_n &=& \big(\big( curl (\phi_n \mathbf{e}_3),\,\phi_n\mathbf{e}_3\big)\big) , 
    \ \ \ \phi_n=\phi_n(r ,\theta ), \\
   -\triangle \phi_n &=& \mu_n \phi_n,\ \phi_n|_{\partial\Gamma}=0, \ \   \mbox{and }\\
       curl\mathbf{\Phi}_n &=&  \big(\big( curl (\phi_n\mathbf{e}_3),\,
         \mu_n\phi_n\mathbf{e}_3\big)\big) .
 \end{eqnarray*}
Here $\big(\big( \mathbf{a},\, b\mathbf{e}_3\big)\big)$
 denotes a 3-component vector whose horizontal projection 
is $\mathbf{a}$
and vertical projection is $b\mathbf{e}_3$.
\end{rem}

Let us  explicit elementary swirling wave flows which are \textit{exact} solutions
to (\ref{nsch8_a}) and (\ref{swiring}):

\begin{lem}\label{lemma2.3}
  For every $n=(n_1,n_2,n_3)$, the following quasiperiodic ($\mathbf{T}^2$ in time)
solenoidal fields are exact solution of the full 3D nonlinear Euler equations
(\ref{nsch8_a}):
\begin{equation}\label{euler_sol}
   \mathbf{V}(t,y)=\frac{\Omega}{2}\mathbf{e}_3\times y + \exp (\frac{\Omega}{2}\mathbf{J}t)
    \mathbf{\Phi}_n(\exp (-\frac{\Omega}{2}\mathbf{J}t)y) 
         \exp (\pm i\frac{n_3}{|\lambda_n|}\alpha\Omega t),
\end{equation}
$n_3$ is the vertical wave number
of $\mathbf{\Phi}_n$ and $\exp (\frac{\Omega}{2}\mathbf{J}t)$ the unitary group
of rigid body rotations:
\begin{equation}
\label{defJ} \mathbf{J}=\left(
\begin{array}{ccc}
0 & -1 & 0 \\
1 & 0 & 0 \\
0 & 0 & 0
\end{array}
\right), \mbox{   }  e^{\Omega{\mathbf{J}}t/2}= \left(
\begin{array}{ccc}
\cos(\frac{\Omega t}{2}) & - \sin(\frac{\Omega t}{2}) & 0 \\
 \sin(\frac{\Omega t}{2}) & \cos(\frac{\Omega t}{2}) & 0 \\
0 & 0 & 1
\end{array}
\right).
\end{equation}
\end{lem}

\begin{rem} These fields are exact quasiperiodic, nonaxisymmetric swirling flow solutions of
the 3D Euler equations.   For $n_3\neq 0$, their second components
\begin{equation}\label{second_compo}
  \tilde{\mathbf{V}}(t,y) = \exp (\frac{\Omega}{2}\mathbf{J}t)\mathbf{\Phi}_n
    (\exp (-\frac{\Omega}{2}\mathbf{J}t) y) \exp (\pm \frac{in_3}{|\lambda_n|}\alpha\Omega t)
\end{equation}
are Beltrami flows ($curl \tilde{\mathbf{V}}\times \tilde{\mathbf{V}}\equiv 0$) exact
solutions of (\ref{swiring}) with $\tilde{\mathbf{V}}(t=0, y) = \mathbf{\Phi}_n(y)$.
\end{rem}

$\tilde{\mathbf{V}}(t,y)$ in Eq. (\ref{second_compo}) are dispersive waves with
frequencies $\frac{\Omega}{2}$ and $\frac{n_3\alpha}{|\lambda_n|}\Omega$, where
$\alpha =\frac{2\pi}{h}$. Moreover, each $\tilde{\mathbf{V}}(t,y)$ is a \textit{traveling wave}
along the cylinder's axis, since it contains the factor 
\[ \exp \left( i\alpha n_3(\pm z \pm \frac{\Omega}{|\lambda_n|}t)\right) .
\]
Note that $n_3$ large corresponds to small axial (vertical) scales, albeit  $0\leq \alpha |n_3/\lambda_n|\leq 1$.

 \textit{Proof of Lemma \ref{lemma2.3}}. Through the canonical rigid body transformation for both the field $\mathbf{V}(t,y)$ and the space coordinates $y=(y_1,y_2,y_3)$:
\begin{equation}
\label{transffund}
{\mathbf{V}}(t,y)=e^{+\Omega {\mathbf{J}}t/2}{\mathbf{U}}(t,e^{-\Omega {%
\mathbf{J}}t/2}y)+\frac{\Omega }{2} {\mathbf{J}} y, \mbox{   }
x=e^{-\Omega {\mathbf{J}}t/2} y,
\end{equation}
the 3D Euler equations (\ref{nsch8_a}), (\ref{nsch8_a/ic}) transform into:
\begin{eqnarray}
\label{nsch8/again}
\partial_t \mbox{\bf{U}} + (curl \mbox{\bf{U}} + \Omega \mathbf{e}_3) \times \mbox{\bf{U}}
 =
- \nabla \left( p - \frac{\Omega^2}{4} (|x_1|^2+|x_2|^2) +
\frac{|\mbox{\bf{U}}|^2}{2}\right), && \\
\label{nsch8/again/ic} \nabla \cdot \mbox{\bf{U}}=0, \ \ \ 
\mbox{\bf{U}}(t,x)|_{t=0} =
\mbox{\bf{U}}(0)=\tilde{\mbox{\bf{V}}}_0(x), \hspace{.7in}&&
\end{eqnarray}

For Beltrami flows such that $curl \mathbf{U}\times \mathbf{U}\equiv 0$,
these Euler equations (\ref{nsch8/again})-(\ref{nsch8/again/ic}) in a rotating frame reduce to:
\[ \partial_t\mathbf{U}+\Omega \mathbf{e}_3\times \mathbf{U} + \nabla \pi = 0, \ \nabla\cdot\mathbf{U}=0,
\]
which are identical to the Poincar\'{e}-Sobolev nonlocal wave equations in the
cylinder \cite{M-N-B-G}, \cite{Poi}, \cite{Sob}, \cite{Ar-Khe}:
\begin{eqnarray}
   \partial_t \boldsymbol{\Psi} + \Omega \mathbf{e}_3\times \mathbf{\Psi}+ 
                \nabla \pi = 0,\ \nabla\cdot\mathbf{\Psi}=0, \hspace{.3in}& & \label{PSwave}\\
   \frac{\partial^2}{\partial t^2} curl^2\mathbf{\Psi}
          -\mathbf{\Omega}^2 \frac{\partial^2}{\partial x_3^2}\mathbf{\Psi} =0, \ 
     \mathbf{\Psi}\cdot \mathbf{N}|_{\partial D}=0. && 
\end{eqnarray}
It suffices to verify that the Beltrami flows $\mathbf{\Psi}_n(t,x) =
   \mathbf{\Phi}_n(x)\exp \left( \pm i\frac{\alpha n_3}{|\lambda_n|}\Omega t\right)$,
where $\mathbf{\Phi}_n^\pm (x)$ and $\pm |\lambda_n|$ are curl eigenfunctions 
and eigenvalues, are exact solutions to the Poincar\'{e}-Sobolev wave equation,
in such a rotating frame of reference.
\begin{rem} The frequency spectrum of the Poincar\'{e} vorticity waves
(solutions to {\em (\ref{PSwave})}) is  exactly $\pm i\frac{\alpha n_3}{|\lambda_n|}\Omega,\ 
  n=(n_1,n_2,n_3)$ indexing the spectrum of curl. 
  Note that $n_3=0$ (zero frequency of rotating waves) corresponds to 2-Dimensional,
3-Components solenoidal vector fields.
\end{rem}

We now transform the Cauchy problem for the 3D Euler equations
 (\ref{nsch8_a})-(\ref{nsch8_a/ic}) into an infinite dimensional 
nonlinear dynamical system by expanding $\mathbf{V}(t,y)$ along the swirling wave
flows (\ref{euler_sol})-(\ref{second_compo}):
{\small
\begin{subequations}
    \begin{align}
   \mathbf{V}(t,y)=\frac{\Omega}{2} \mathbf{e}_3\times y  \hspace{3.4in}
   \\+
      \exp \left( \frac{\Omega}{2}\mathbf{J}t \right)
      \left\{ \sum_{n=(n_1,n_2,n_3)}\mathbf{u}_n(t)
              \exp\left( \pm i\frac{\alpha n_3}{|\lambda_n|}\Omega t\right)
               \mathbf{\Phi}_n\left( \exp \left( -\frac{\Omega}{2}\mathbf{J}t\right) y\right)
       \right\} \hspace{.2in} \\
    \mathbf{V}(t=0,y) = \frac{\Omega}{2}\mathbf{e}_3\times y+\tilde{\mathbf{V}}_0(y)
            \hspace{2in}  \\
    \tilde{\mathbf{V}}_0(y) =
       \sum_{n=(n_1,n_2,n_3)}\mathbf{u}_n(0)\mathbf{\Phi}_n(y), \hspace{1.9in} 
             \label{eq2.23c}
    \end{align}
\end{subequations}
}
where $\mathbf{\Phi}_n$ denotes the curl eigenfunctions of Proposition 2.1
if $n_3\neq 0$, and $\mathbf{\Phi}_n = \big(\big( curl (\phi_n\mathbf{\mathbf{e}_3}),
  \,\phi_n\mathbf{e}_3\big)\big)$ if
$n_3=0$ (2D case, Remark 2.2).
\medskip

As we focus on the case where helicity is conserved for (\ref{nsch8_a})-(\ref{nsch8_a/ic}),
we consider the class of initial data $\tilde{\mathbf{V}}_0$ such that
\cite{M-N-B-G}:
\[ curl\tilde{\mathbf{V}}_0\cdot\mathbf{N} = 0 \ \mbox{on}\ \Gamma ,\]
where $\Gamma$ is the lateral boundary of the cylinder.

\medskip

The infinite dimensional dynamical system is then equivalent to the 3D Euler
equations (\ref{nsch8_a})-(\ref{nsch8_a/ic}) in the cylinder,
with $n=(n_1,n_2,n_3)$ ranging over the whole spectrum of curl, e.g.:
\begin{equation}
\begin{aligned}
\frac{d\mathbf{u}_n}{dt}=- \sum_{\substack{k,m\\k_3+m_3=n_3\\ k_2+m_2=n_2}}
  \exp \left( i \left( \pm \frac{n_3}{|\lambda_n|}\pm\frac{k_3}{|\lambda_k|}
     \pm\frac{m_3}{|\lambda_m|}\right) \alpha \Omega t\right)
\\
\times < curl\mathbf{\Phi}_k\times\mathbf{\Phi}_m, \mathbf{\Phi}_n>
\mathbf{u}_k(t)\mathbf{u}_m(t)\label{eq2.29} 
\end{aligned}
\end{equation}

\medskip

\noindent here 
\begin{subequations}
  \begin{align}
   curl \mathbf{\Phi}_k^{\pm} = \pm \lambda_k\mathbf{\Phi}_k^{\pm}
  \ \mbox{ if }\ k_3\neq 0, \nonumber \\
  curl \mathbf{\Phi}_k= \big(\big( curl (\phi_k\mathbf{e}_3),
   \,\mu_k\phi_k\mathbf{e}_3\big)\big) \ \mbox{ if }\ k_3=0\nonumber
  \end{align}
\end{subequations}
 (2D, 3-components, Remark 2.2), similarly for $m_3=0$ and $n_3=0$.
The inner product $<\ ,\ >$ denotes the $L^2$ complex-valued  inner product in $D$.

This is an infinite dimensional system of coupled equations  with quadratic
nonlinearities, which conserve both the energy
\[ E(t) = \sum_n |\mathbf{u}_n(t)|^2
\]
 and the helicity
\[ Hel (t) = \sum_{n}\pm |\lambda_n|\,|\mathbf{u}_n^{\pm}(t)|^2.
\]
The quadratic nonlinearities split into \textit{resonant terms}
where the exponential oscillating phase factor in (\ref{eq2.29}) reduces
to unity and fast oscillating non-resonant terms ($\Omega >> 1$).
The resonant set $K$ is defined
in terms of vertical wavenumbers $k_3, m_3, n_3$ and eigenvalues
$\pm \lambda_k$, $\pm \lambda_m$, $\pm \lambda_n$ of $curl$:
\begin{eqnarray}
\label{rescc} K=\{ \pm \frac{k_3}{\lambda_k} \pm
\frac{m_3}{\lambda_m} \pm \frac{n_3}{\lambda_n}=0, n_3=k_3+m_3,
n_2=k_2+m_2 \}.
\end{eqnarray}
Here $k_2, m_2, n_2$ are azimuthal wavenumbers.

We shall call the ``resonant Euler equations'' the following 
$\infty$-dimensional dynamical system restricted to
$(k,m,n)\in K$:
\begin{subequations}\label{resonant_euler}
  \begin{align}
  \frac{d\mathbf{u}_n}{dt}+\sum_{(k,m,n)\in K}
    < curl\mathbf{\Phi}_k\times \mathbf{\Phi}_m, \mathbf{\Phi}_n > 
         \mathbf{u}_k\mathbf{u}_m = 0,  \label{resonant_euler_a}\\
    \mathbf{u}_n(0) \equiv <\tilde{\mathbf{V}}_0, \mathbf{\Phi}_n >, \hspace{1.4in}
            \label{eq2.26b}
   \end{align}
\end{subequations}
\medskip

\noindent here $curl \mathbf{\Phi}_k^{\pm} = \pm \lambda_k\mathbf{\Phi}_k^{\pm}$
if $k_3\neq 0$, 
$curl\mathbf{\Phi}_k= \big(\big( curl (\phi_k\mathbf{e}_3),
   \,\mu_k\phi_k\mathbf{e}_3\big)\big)$ if $k_3=0$; 
similarly for $m_3=0$ and $n_3=0$ (2D components, Remark 2.2).
If there are no terms in (\ref{resonant_euler_a}) satisfying the resonance conditions, then 
there will be some modes for which
$\frac{d\mathbf{u}_j}{dt}=0$.

\begin{lem}
 The resonant 3D Euler equations (\ref{resonant_euler}) 
 conserve both energy $E(t)$ and helicity $Hel(t)$. The energy and 
helicity are identical to that of the
  full exact 3D Euler equations (\ref{nsch8_a})-(\ref{nsch8_a/ic}).
\end{lem}

The set of resonances $K$ is studied in depth in \cite{M-N-B-G}.
To summarize, $K$ splits into:
\begin{itemize}
  \item[(i\,)] $0$-wave resonances, with $n_3=k_3=m_3=0$; the corresponding resonant equations
 are identical to the 2-Dimensional, 3-Components Euler equations, with I.C.
 \[ \frac{1}{h}\int_0^h\tilde{\mathbf{V}}_0(y_1,y_2,y_3)\,dy_3. 
\]
  \item[(ii)] Two-Wave resonances, with $k_3m_3n_3=0$, but two of them are not null;
  the corresponding resonant equations (called ``catalytic equations'') are 
   proven to possess an infinite, countable set of new conservation laws 
   \cite{M-N-B-G}.
  \item[(iii)] Strict three-wave resonances for a subset $K^*\subset K$.
\end{itemize}

\begin{defn}
   The set $K^*$ of strict 3 wave resonances is:
   \begin{equation}\label{strict_resonances}
      K^*=\left\{ \pm\frac{k_3}{\lambda_k}\pm\frac{m_3}{\lambda_m}\pm\frac{n_3}{\lambda_n}
       =0, k_3m_3n_3\neq 0, n_3=k_3+m_3, n_2=k_2+m_2 \right\} .
   \end{equation}
\end{defn}
Note that $K^*$ is parameterized by $h/R$, since $\alpha=\frac{2\pi}{h}$
parameterizes the eigenvalues $\lambda_n, \lambda_k, \lambda_m$ of the curl operator.

\begin{prop}
 There exist a countable, non-empty set of parameters $\frac{h}{R}$ for which $K^*\neq \varnothing$.
\end{prop}

\textit{Proof.} The technical details, together with a more precise statement,
 are postponed to the proof of Lemma \ref{lemma5.7}. Concrete examples of
resonant axisymmetric and helical waves 
are discussed in \cite{Mah} ( cf. Figure 2 in the article).

\begin{cor}
  Let $\int_0^h \tilde{\mathbf{V}}_0(y_1,y_2,y_3)\,dy_3 =0$, i.e. zero vertical
  mean for the I.C. $\tilde{\mathbf{V}}_0(y)$ in (\ref{nsch8_a/ic}), (\ref{eq2.8}),
  (\ref{eq2.23c}) and (\ref{eq2.26b}). Then the resonant 3D Euler equations are 
 invariant on $K^*$:
 \begin{subequations}\label{resonant_euler_again}
   \begin{align}
     \frac{d\mathbf{u}_n}{dt}+\sum_{(k,m,n)\in K^*}\lambda_k
     <\mathbf{\Phi}_k\times \mathbf{\Phi}_m, \mathbf{\Phi}_n>\mathbf{u}_k\mathbf{u}_m=0 ,\ 
         k_3m_3n_3\neq 0, \\
     \mathbf{u}_n(0) = < \tilde{\mathbf{V}}_0, \mathbf{\Phi}_n> \hspace{1in}
   \end{align}
 \end{subequations}
(where $\tilde{\mathbf{V}}_0$ has spectrum restricted to $n_3\neq 0$).
\end{cor}

\textit{Proof.} This is an immediate corollary of
the ``operator splitting'' Theorem 3.2 in \cite{M-N-B-G}.
$\ \ \blacksquare$

We shall call the above dynamical systems the ``\textit{strictly resonant Euler system}''.
This is an $\infty$-dimensional Riccati system which conserves Energy and Helicity.
It corresponds to \textit{nonlinear interactions depleted} on $K^*$.

How do dynamics of the resonant Euler equations (\ref{resonant_euler}) or 
(\ref{resonant_euler_again}) approximate \textit{exact} solutions of the Cauchy
problem for the full Euler equations in strong norms?
This is answered by the following theorem, proven in Section 4 of 
\cite{M-N-B-G}:
\begin{thm}\label{theorem2.10} Consider the initial value problem
\[ \mathbf{V}(t=0,y) = \frac{\Omega}{2} \mathbf{e}_3\times y +\tilde{\mathbf{V}}_0(y),\ 
   \tilde{\mathbf{V}}_0\in \mathbf{H}_{\sigma}^s,\ s > 7/2
\]
for the full 3D Euler equations, with $||\tilde{\mathbf{V}}_0||_{\mathbf{H}_{\sigma}^s}\leq M_s^0$
and $curl \tilde{\mathbf{V}}_0\cdot \mathbf{N} = 0$ on~$\Gamma$.
\begin{itemize}
 \item Let $\mathbf{V}(t,y)=\frac{\Omega}{2}\mathbf{e}_3\times y + \tilde{\mathbf{V}}(t,y)$
     denote the solution to the \textit{exact} Euler equations.
 \item Let $\mathbf{w}(t,x)$ denote the solution to the \textit{resonant} 3D Euler equations
with Initial Condition $\mathbf{w}(0, x)\equiv \mathbf{w}(0,y)=\tilde{\mathbf{V}}_0(y)$.
 \item Let $||\mathbf{w}(t,y)||_{\mathbf{H}_{\sigma}^s}\leq M_s(T_M,M_s^0)$ on $0\leq t\leq T_M,\ s>7/2$.
\end{itemize}
   Then, $\forall \epsilon > 0 ,\ \exists \ \Omega^*(T_M, M_s^0, \epsilon )$
   such that, $\forall \Omega \geq \Omega^*$:
   \[ \Big|\Big|\tilde{\mathbf{V}}(t,y) - 
       \exp \left( \frac{\mathbf{J}\Omega t}{2}\right)
              \left\{ \sum_n \mathbf{u}_n(t) e^{-i\frac{n_3}{\lambda_n}\alpha\Omega t}
                \mathbf{\Phi}_n (e^{-\frac{\mathbf{J}\Omega t}{2}}y)\right\}  \Big|\Big|_{H^\beta}
   \leq \epsilon 
\]
 on $0\leq t\leq T_M,\ 
   \forall \beta\geq 1,\ \beta \leq s-2$. Here $||\cdot ||_{H^\beta}$ is defined
 in (\ref{Hs}).
\end{thm}
The 3D Euler flow preserves the condition 
$curl \tilde{\mathbf{V}}_0\cdot \mathbf{N} = 0$ on~$\Gamma$,
that is $curl \mathbf{V}(t,y)\cdot\mathbf{N}=0$ on $\Gamma$, for every $t\geq 0$
\cite{M-N-B-G}.
The proof of this ``error-shadowing'' theorem is delicate, beyond the usual
Gronwall differential inequalities and involves estimates of oscillating integrals
of almost periodic functions of time with values in Banach spaces.
Its importance lies in that solutions of the resonant Euler equations (\ref{resonant_euler})
and/or (\ref{resonant_euler_again}) are uniformly close  in strong norms
to those of the exact Euler equations (\ref{nsch8_a})-(\ref{nsch8_a/ic}), on any time interval of
existence of smooth solutions of the \textit{resonant} system. The infinite dimensional
Riccati systems (\ref{resonant_euler})
and (\ref{resonant_euler_again}) are not just hydrodynamic models, but \textit{exact asymptotic
limit systems} for $\Omega \gg 1$. This is  in contrast to all previous
literature on conservative 3D hydrodynamic models, such as in \cite{G-D-O}.

%-----------------------------------------------------
  \section{Strictly resonant Euler systems: the SO(3) case}
%-----------------------------------------------------

We investigate the structure and the dynamics of the ``strictly resonant
Euler systems'' (\ref{resonant_euler_again}). Recall that
the set of 3-wave resonances is:
\begin{equation}
\begin{aligned}
 K^*=\Big\{
     (k,m,n): \pm \frac{k_3}{\lambda_k}\pm \frac{m_3}{\lambda_m}\pm \frac{n_3}{\lambda_n}=0,
     k_3m_3n_3\neq 0, 
\\
n_3=k_3+m_3, n_2=k_2+m_2 \Big\} .
\end{aligned}
\end{equation}
From the symmetries of the curl eigenfunctions $\mathbf{\Phi}_n$ and eigenvalues 
$\lambda_n$ in the cylinder,  the following identities hold under the transformation 
$n_2\rightarrow -n_2$, $n_3\rightarrow -n_3$
\begin{equation}
\begin{aligned}
\mathbf{\Phi}(n_1, -n_2, -n_3) = \mathbf{\Phi}^*(n_1, n_2, n_3)\,,
\\
\lambda (n_1, -n_2, -n_3)=\lambda (n_1, n_2, n_3)\,.
\end{aligned}
\end{equation}
where $^*$ designates the complex conjugate (see Section 3, \cite{M-N-B-G}
for details). \ The eigenfunctions\ \ $\mathbf{\Phi}(n_1,n_2,n_3)$ \ involve
the radial functions\  \ 
$J_{n_2}(\beta (n_1,n_2,\alpha n_3)r)$ and  $J_{n_2}^\prime (\beta (n_1,n_2,\alpha n_3)r)$,
with 
\[ \lambda^2(n_1,n_2,n_3) = \beta^2(n_1,n_2,\alpha n_3)+\alpha^2 n_3^2;
\]
$\beta (n_1,n_2,\alpha n_3)$ are discrete, countable roots of
equation (3.30) in \cite{M-N-B-G}, obtained via an equivalent Sturm-Liouville radial
problem. Since the curl eigenfunctions are even in $r\rightarrow -r,\
n_1\rightarrow -n_1$, we will extend the indices $n_1=1,2,...,+\infty$
to $-n_1=-1,-2,...$ with the above radial symmetry in mind.
\begin{cor}\label{cor3.1}
The 3-wave resonance set $K^*$ is invariant under the symmetries $\sigma_j,\ 
j=0,1,2,3$, where 
$$
\begin{aligned}
\sigma_0(n_1,n_2,n_3)&=(n_1,n_2,n_3),
\\
\sigma_1(n_1,n_2,n_3)&= (-n_1,n_2,n_3),
\\
\sigma_2(n_1,n_2,n_3)&= (n_1,-n_2,n_3)
\\
\sigma_3(n_1,n_2,n_3)&=(n_1,n_2,-n_3)\,.
\end{aligned}
$$
\end{cor}

\begin{rem}\label{rem3.2}
  For $0<i\leq 3,\ 0<j\leq 3,\ 0<l\leq 3$  $\sigma_j^2=Id,\ \sigma_i\sigma_j=-\sigma_l$
  if $ i\neq j$ and $\sigma_i\sigma_j\sigma_l=-Id$, for $i\neq j\neq l$.
  The $\sigma_j$ do preserve the convolution conditions in $K^*$.
\end{rem}

We choose an $\alpha$ for which the set $K^*$ is not empty. We further take the
hypothesis of a single triple wave resonance $(k,m,n)$, modulo the symmetries $\sigma_j$:
\begin{hyp}\label{hypothesis5.3}
  $K^*$ is such that there exists a single triple wave number resonance
$(n,k,m)$, modulo the symmetries $\sigma_j,\ j=1,2,3$ and 
$\sigma_j(k)\neq k,\ \sigma_j(m)\neq m,\ \sigma_j(n)\neq n$ for
$j=2$ and $j=3$.
\end{hyp}
Under the above hypothesis, one can demonstrate that the strictly resonant
Euler system splits into three uncoupled systems in $\mathbf{C}^3$:
\begin{thm} 
 Under hypothesis \ref{hypothesis5.3}, the resonant Euler system reduces to three
uncoupled rigid body systems in $\mathbf{C}^3$:
\begin{subequations}\label{reduced_rigid}
   \begin{align}
       \frac{dU_n}{dt}+i (\lambda_k-\lambda_m)C_{kmn}U_kU_m = 0 \\ 
       \frac{dU_k}{dt}-i (\lambda_m-\lambda_n)C_{kmn}U_nU_m^* = 0 \\
        \frac{dU_m}{dt}-i (\lambda_n-\lambda_k)C_{kmn}U_nU_k^* = 0 
   \end{align}
\end{subequations}
where $C_{kmn}=i <\mathbf{\Phi}_k\times\mathbf{\Phi}_m, \mathbf{\Phi}_n^*>,\
  C_{kmn}$ real and the other two \textit{uncoupled} systems obtained with the symmetries
$\sigma_2(k,m,n)$ and $\sigma_3(k,m,n)$. The energy and the helicity of each subsystem
are conserved:
\begin{eqnarray*}
  \frac{d}{dt}(U_kU_k^*+U_mU_m^*+U_nU_n^*) = 0, \hspace{.5in}&&\\
  \frac{d}{dt}(\lambda_k U_kU_k^*+\lambda_m U_mU_m^*+\lambda_n U_nU_n^*) = 0. &&
\end{eqnarray*}
\end{thm}

\textit{Proof.} It follows from $U_{-k}=U_k^*,\ \lambda (-k)=\lambda (+k)$,
 similarly for $m$ and $n$; and in a very essential way from the antisymmetry of
$<\mathbf{\Phi}_k\times \mathbf{\Phi}_m, \mathbf{\Phi}_n^*>$, together with
$curl \mathbf{\Phi}_k=\lambda_k\mathbf{\Phi}_k$. That $C_{kmn}$ is real follows
from the eigenfunctions explicited in Section 3 of \cite{M-N-B-G}.$\ \ \blacksquare$

\begin{rem}
  This deep structure, i.e. $SO(3;\mathbf{C})$
 rigid body systems in $\mathbf{C}^3$ is a direct 
consequence of the Lam\'{e} form of the full 3D Euler equations, cf.
Eqs. (\ref{Lame_form}) and (\ref{swiring}), and the nonlinearity $curl \mathbf{V}\times\mathbf{V}$.
\end{rem}

The system (\ref{reduced_rigid}) is equivariant with respect to the symmetry
operators
\[ (z_1,z_2,z_3)\rightarrow (z_1^*,z_2^*,z_3^*),\
   (z_1,z_2,z_3)\rightarrow \left( \exp (i\chi_1)z_1, \exp (i\chi_2)z_2, \exp (i\chi_3)z_3
\right),
\]
  provided $\chi_1=\chi_2+\chi_3$.
It admits other integrals known as the  Manley-Rowe relations (see, for instance 
\cite{We-Wil}). It differs from the usual 3-wave resonance systems investigated 
in the literature, such as in \cite{Zak-Man1}, \cite{Zak-Man2},
 \cite{Gu-Ma} in that
\begin{itemize}
  \item[(1)] helicity is conserved,
  \item[(2)] dynamics of these resonant systems rigorously ``shadow'' those 
    of  the exact 3D Euler equations, see  Theorem \ref{theorem2.10}.
\end{itemize}

Real forms of the system (\ref{reduced_rigid}) are found in Gledzer et al. \cite{G-D-O},
corresponding to the exact invariant manifold $U_k\in i\mathbf{R},\ 
U_m\in\mathbf{R}, \ U_n\in\mathbf{R}$, albeit without any rigorous asymptotic
justification. The $\mathbf{C}^3$ systems (\ref{reduced_rigid})  with
\textit{helicity} conservation laws are not discussed in
\cite{G-D-O}.

The only nontrivial Manley-Rowe conservation laws for the resonant system
(\ref{reduced_rigid}), rigid body $SO(3;\mathbf{C})$, which are
 independent from energy and helicity, are:
\[ \frac{d}{dt}\left( r_kr_mr_n \sin (\theta_n-\theta_k-\theta_m)\right) = 0,
\]
where $U_j=r_j \exp (i\theta_j),\ j=k,m,n,$ and
\begin{eqnarray*}
   \mathcal{E}_1 = (\lambda_k-\lambda_m) r_n^2-(\lambda_m-\lambda_n) r_k^2, && \\
   \mathcal{E}_2 = (\lambda_m-\lambda_n) r_k^2-(\lambda_n-\lambda_k) r_m^2. && 
\end{eqnarray*}

The resonant system (\ref{reduced_rigid}) is well known to possess hyperbolic equilibria
and heteroclinic/homoclinic orbits on the energy surface.
We are interested in rigorously proving arbitrary large bursts of enstrophy and higher
norms on arbitrarily small time intervals,
for properly chosen $h/R$. To simplify the presentation, we establish
the results for the simpler invariant manifold
$U_k\in i\mathbf{R}$, and $U_m, U_n \in \mathbf{R}$.

Rescale time as:
\[ t \rightarrow t/C_{kmn}.
\]

Start from the system
\begin{equation}
\label{SysStrt}
\begin{aligned}
\dot{U}_n+i(\l_k-\l_m)U_kU_m&=0
\\
\dot{U}_k-i(\l_m-\l_n)U_nU_m^*&=0
\\
\dot{U}_m-i(\l_n-\l_k)U_nU_k^*&=0
\end{aligned}
\end{equation}
Assume that $U_k\in i\bR$ and that $U_m,U_n\in\bR$: set $p=iU_k$, $q=U_m$ and $r=U_n$, as well as
$\l_k=\l$, $\l_m=\mu$ and $\l_n=\nu$: then
\begin{equation}
\label{SysRl}
\begin{aligned}
\dot{p}+(\mu-\nu)qr&=0
\\
\dot{q}+(\nu-\l)rp&=0
\\
\dot{r}+(\l-\mu)pq&=0
\end{aligned}
\end{equation}
This system admits two first integrals:
\begin{equation}
\label{1Int}
\begin{aligned}
E&=p^2+q^2+r^2\qquad\qquad&&\hbox{(energy)}
\\
H&=\l p^2+\mu q^2+\nu r^2&&\hbox{(helicity)}
\end{aligned}
\end{equation}

System (\ref{SysRl}) is exactly the $SO(3, \mathbf{R})$ rigid body dynamics Euler equations,
with inertia momenta $I_j=\frac{1}{|\lambda_j|},\ j=k,m,n$ \cite{Ar1}.
\begin{lem} {\em (\cite{Ar1}, \cite{G-D-O})} 
  With the ordering $\lambda_k > \lambda_m > \lambda_n$, i.e. $\lambda > \mu > \nu$, 
the equilibria $(0, \pm 1, 0)$ are hyperbolic saddles on the unit energy sphere, and
the equilibria $(\pm 1, 0, 0),\ (0,0,\pm 1)$ are
centers. There exist equivariant families of 
 heteroclinic connections
between $(0, +1, 0)$ and $(0, -1, 0)$. Each pair of such connections correspond to
equivariant homoclinic cycles at $(0,1,0)$ and $(0, -1, 0)$.
\end{lem}

We investigate bursting dynamics along orbits with large periods, with initial 
conditions close to the hyperbolic point $(0, E(0), 0)$ on the energy
sphere $E$. We choose resonant triads such that $\lambda_k > 0, \ \lambda_n<0,
\lambda_k\sim |\lambda_n|,\ |\lambda_m|\ll \lambda_k,$ equivalently:
\begin{equation}
   \lambda > \mu > \nu,\ \lambda \nu <0,\ |\mu|\ll \lambda \mbox{ and } \lambda\sim |\nu|.
\end{equation}

\begin{lem}\label{lemma5.7}
  There exist $h/R$ with $K^*\neq \varnothing$, such that
\[ \lambda_k >\lambda_m >\lambda_n,\ \lambda_k\lambda_n <0,\ |\lambda_m|\ll \lambda_k
  \mbox{ and } \lambda_k\sim |\lambda_n|.
\]
\end{lem}
\begin{rem}
  Together with the polarity $\pm$ of the curl eigenvalues,
these are 3-wave resonances where two of the eigenvalues are much
larger in moduli than the third one. In the limit $|k|,\ |m|, \ |n|\gg 1,
\lambda_k\sim\pm |k|,\ \lambda_m\sim \pm |m|,\ \lambda_n\sim \pm |n|,$
the eigenfunctions $\mathbf{\Phi}$ have leading asymptotic terms which involve
 cosines and sines periodic in $r$, cf. Section 3 {\em \cite{M-N-B-G}}.
In the strictly resonant equations (\ref{resonant_euler_again}),
the summation over the quadratic terms becomes an asymptotic convolution in
$n_1=k_1+n_1$. The resonant three waves in Lemma \ref{lemma5.7} are
equivalent to Fourier triads $k+m=n,$ with $|k|\sim |n|$ and $|m|\ll |k|, \ |n|$,
in periodic lattices. In the physics of spectral theory of turbulence 
{\em \cite{Fri}, \cite{Les}}, 
these are exactly the triads responsible from transfer of energy between large scales
and small scales. These are  the triads which
have hampered  mathematical efforts at proving the global regularity
of the Cauchy problem for 3D Navier-Stokes equations in periodic lattices
{\em \cite{Fe}}.
\end{rem}

\textit{Proof of Lemma \ref{lemma5.7}} (\cite{M-N-B-G}) The transcendental dispersion law
for 3-waves in $K^*$ for cylindrical domains, is a polynomial
of degree four in $\vartheta_3=1/h^2$:
\begin{equation}\label{trancendental}
  \tilde{P}(\vartheta_3)=\tilde{P}_4\vartheta_3^4+\tilde{P}_3\vartheta_3^3
    +\tilde{P}_2\vartheta_3^2+\tilde{P}_1\vartheta_3 + \tilde{P}_0 = 0,
\end{equation}
with $n_2=k_2+m_2$ and $n_3=k_3+m_3$. 

\nopagebreak[1]
Then with $h_k=\frac{\beta^2(k_1,k_2,\alpha k_3)}{k_3^2},\ 
  h_m=\frac{\beta^2(m_1,m_2,\alpha m_3)}{m_3^2},\
  h_n=\frac{\beta^2(n_1,n_2,\alpha n_3)}{n_3^2}$,  
 cf. the radial Sturm-Liouville problem
in Section 3, \cite{M-N-B-G},  the coefficients of $\tilde{P}(\vartheta_3)$
are given by:
\begin{eqnarray*}
   \tilde{P}_4 &=& -3, \\
   \tilde{P}_3 &=& -4 (h_k+h_m+h_n), \\
   \tilde{P}_2 &=& -6 (h_kh_m+h_kh_n+h_mh_n), \\
   \tilde{P}_1 &=& -12 h_kh_mh_n, \\
   \tilde{P}_0 &=& h_m^2h_n^2+h_k^2h_n^2+h_m^2h_k^2
                -2(h_kh_mh_n^2+h_kh_nh_m^2+h_mh_nh_k^2). 
\end{eqnarray*}
Similar formulas for the periodic lattice domain were first derived in \cite{B-M-N2},
\cite{B-M-N3},
\cite{B-M-N4}. 
In cylindrical domains the resonance condition for $K^*$ is identical to
\[  \pm \frac{1}{\sqrt{\vartheta_3 +h_k}}
    \pm \frac{1}{\sqrt{\vartheta_3 +h_m}}
    \pm \frac{1}{\sqrt{\vartheta_3 +h_n}} = 0,
\]
with $\vartheta_3=\frac{1}{h^2},\ h_k=\beta^2(k)/k_3^2,\ 
  h_m=\beta^2(m)/m_3^2,\ h_n=\beta^2(n)/n_3^2$; 
Eq. (\ref{trancendental}) is the equivalent rational form.

From the asymptotic formula (3.44) in \cite{M-N-B-G}, 
for large $\beta$:
\begin{equation}
 \beta (n_1,n_2,n_3)\sim n_1\pi +n_2\frac{\pi}{2}+\frac{\pi}{4}+\psi ,
\end{equation}
where $\psi = 0$ if $\lim \frac{m_2}{m_3}\frac{h}{2\pi}=0$ (e.g.
$h$ fixed, $m_2/m_3\rightarrow 0$) and $\psi = \pm \frac{\pi}{2}$ if $\lim \frac{m_2}{m_3}\frac{h}{2\pi}=\pm\infty$
(e.g. $\frac{m_2}{m_3}$ fixed, $h\rightarrow\infty$).
The proof is completed by taking leading terms
$\tilde{P}_0+\vartheta_3\tilde{P}_1$ in (\ref{trancendental}), 
$\vartheta_3=\frac{1}{h^2}\ll 1$,
and $m_2=0,\ k_2=\mathcal{O}(1),\ n_2=\mathcal{O}(1). \ \ \ \blacksquare$
\medskip

We now state a theorem for bursting of the $\mathbf{H}^3$ norm in arbitrarily 
small times, for initial data close to the hyperbolic point $(0, E(0), 0)$:
\begin{thm}\label{thm3.9} (Bursting dynamics in $\mathbf{H}^3$). 
  Let $\lambda > \mu > \nu ,\ \lambda\nu <0,\ |\mu|\ll \lambda$ and $\lambda\sim |\nu|$.
Let $W(t)= \lambda^6 p(t)^2+\mu^6 q(t)^2+\nu^6 r(t)^2$ the $\mathbf{H}^3$-norm squared of
an orbit of (\ref{SysRl}). Choose initial data such that:
$ W(0)=\lambda^6 p(0)^2+\mu^6 q(0)^2$ with $ \lambda^6 p(0)^2\sim \frac{1}{2} W(0)$ and
$\mu^6 q(0)^2\sim\frac{1}{2}W(0)$. Then there exists $t^*>0$, such that
\[ W(t) \geq \frac{1}{4}\left( \frac{\lambda}{\mu}\right)^6 W(0) \]
where $t^*\leq \frac{6}{\sqrt{W(0)}}\mu^2 Ln(\lambda /|\mu|)(\lambda /|\mu|)^{-1}$.
\end{thm}

\begin{rem}
Under the conditions of Lemma \ref{lemma5.7}, $\left(\frac{\lambda}{\mu}\right)^6\gg 1$, 
whereas $\mu^2 (Ln (\lambda /|\mu| ) )(\lambda /|\mu|)^{-1}\ll 1$. Therefore, over a small 
time interval of length $O(\mu^2 (Ln (\lambda /|\mu| ) )(\lambda /|\mu|)^{-1})\ll 1$, the ratio
$||\mathbf{U}(t)||_{H^3}/||\mathbf{U}(0)||_{H^3}$ grows up to a maximal value 
$O\left((\lambda/|\mu|)^3\right)\gg 1$. Since the orbit is periodic, the $\mathbf{H}^3$ 
semi-norm  eventually relaxes to its initial state after some time (this being a 
manifestation of the time-reversibility of the Euler flow on the energy sphere). The 
``shadowing'' theorem \ref{theorem2.10} with $s> 7/2$ ensures that the full, original
3D Euler dynamics, with the same initial conditions, will undergo the same type of 
burst. Notice that, with the definition (\ref{Hs}) of $\|\cdot\|_{H^s}$, one has 
$$
||\Omega \mathbf{e}_3\times y||_{H^3}=||curl^3(\Omega \mathbf{e}_3\times y)||_{L^2}=0\,.
$$
Hence the solid rotation part of the original 3D Euler solution does not contribute to
the ratio $||\mathbf{V}(t)||_{H^3}/||\mathbf{V}(0)||_{H^3}$.
\end{rem}

\begin{thm}\label{thm3.11} (Bursting dynamics of the enstrophy). Under the same conditions for 
the 3-wave resonance, let $\Xi (t) = \lambda^2p(t)^2+\mu^2 q(t)^2+\nu^2 r(t)^2$ the enstrophy.
Choose initial data such that $\Xi (0)=\lambda^2 p(0)^2+\mu^2 q(0)^2+\nu^2 r(0)^2$ with
$\lambda^2 p(0)^2\sim \frac{1}{2} \Xi (0),\ 
\mu^2 q(0)^2\sim \frac{1}{2} \Xi (0)$.
Then there exists $t^{**}>0$, such that
\[ \Xi (t^{**}) \geq \left( \frac{\lambda}{\mu}\right)^2
\]
where $t^{**}\leq \frac{1}{\sqrt{2}}\frac{1}{\sqrt{\Xi (0)}} 
   Ln\left( \lambda /|\mu|\right) \left( \lambda /|\mu|\right)^{-1}.$
\end{thm}

\begin{rem}
It is interesting to compare this mechanism for bursts with earlier results in the
same direction obtained by DiPerna and Lions. Indeed, for each $p\in(1,\infty)$,
each $\delta\in(0,1)$ and each $t>0$, Di Perna and Lions  {\em \cite{DiPe-Li}} 
constructed examples of 2D-3 components solutions to Euler equations such 
that 
$$
||\mathbf{V}(0)||_{W^{1,p}}\leq \epsilon\quad
	\hbox{ while }\quad||\mathbf{V}(t)||_{W^{1,p}}\geq 1/\delta\,.
$$
Their examples essentially correspond to shear flows of the form
$$
\mathbf{V}(t,x_1,x_2)=\left(
\begin{matrix}u(x_2)\\0\\w(x_1-tu(x_2),x_2)\end{matrix}\right)
$$
where $u\in W^{1,p}_{x_2}$ while $w\in W^{1,p}_{x_2}$. Obviously
$$
curl\mathbf{V}(t,x_1,x_2)=\left(\begin{matrix}
(\partial_2-tu'(x_2)\partial_1)w(x_1-tu(x_2),x_2)\\
-\partial_1w(x_1-tu(x_2),x_2)\\
-u'(x_2)\end{matrix}\right)
$$
Thus, all components in $curl\mathbf{V}(t,x_1,x_2)$ belong to $L^p_{loc}$,
except for the term
$$
-tu'(x_2)\partial_1w(x_1-tu(x_2),x_2)\,.
$$
For each $t>0$, this term belongs to $L^p$ for all choices of the functions 
$u\in W^{1,p}_{x_2}$ and $w\in W^{1,p}_{x_1,x_2}$ if and only if $p=\infty$.
Whenever $p<\infty$, DiPerna and Lions construct their examples as some
smooth approximation of the situation above in the strong $W^{1,p}$ topology.

In other words, the DiPerna-Lions construction works only in cases where 
the initial vorticity does not belong to an algebra --- specifically to $L^p$,
which is not an algebra unless $p=\infty$.

The type of burst obtained in our construction above is different: in that
case, the original vorticity belongs to the Sobolev space $H^2$, which is
an algebra in space dimension $3$. Similar phenomena are observed
in all Sobolev spaces $H^\beta$ with $\beta\ge 2$ --- which are also
algebras in space dimension 3. 

In other words, our results complement those of DiPerna-Lions on bursts
in higher order Sobolev spaces, however at the expense of using more 
intricate dynamics.
\end{rem}

We proceed to the proofs of Theorem \ref{thm3.9} and \ref{thm3.11}.
We are interested in the evolution of
\begin{equation}
\label{DfntOm}
\Xi =\l^2 p^2+\mu^2 q^2+\nu^2 r^2\qquad\hbox{(enstrophy)}
\end{equation}
Compute
\begin{equation}
\label{Omdot=}
\dot{\Xi}=-2\left(\l^2(\mu-\nu)+\mu^2(\nu-\l)+\nu^2(\l-\mu)\right)pqr
\end{equation}
then
\begin{equation}
\label{pqrdot}
\dot{(pqr)}=-(\mu-\nu)q^2r^2-(\nu-\l)r^2p^2-(\l-\mu)p^2q^2
\end{equation}
Using the first integrals above, one has
\begin{equation}
\label{pqr>EHOm}
Van\left(\begin{matrix} p^2\\q^2\\r^2\end{matrix}\right)
=
\left(\begin{matrix} E\\H\\ \Xi\end{matrix}\right)
\end{equation}
where $Van$ is the Vandermonde matrix
$$
Van=\left(\begin{matrix}
1&1&1\\
\l     &\mu    &\nu\\
\l^2     &\mu^2    &\nu^2
\end{matrix}\right)
$$
For $\l\not=\mu\not=\nu\not=\l$, this matrix is invertible and
$$
Van^{-1}=\left(\begin{matrix}
\frac{\mu\nu}{(\l-\mu)(\l-\nu)}&\frac{-(\mu+\nu)}{(\l-\mu)(\l-\nu)}&\frac{1}{(\l-\mu)(\l-\nu)}\\
\frac{\nu\l}{(\mu-\nu)(\mu-\l)}&\frac{-(\nu+\l)}{(\mu-\nu)(\mu-\l)}&\frac{1}{(\mu-\nu)(\mu-\l)}\\
\frac{\l\mu}{(\nu-\l)(\nu-\mu)}&\frac{-(\l+\mu)}{(\nu-\l)(\nu-\mu)}&\frac{1}{(\nu-\l)(\nu-\mu)}
\end{matrix}\right)
$$
Hence
\begin{equation}
\label{EHOm>pqr}
\begin{aligned}
p^2&=\frac1{(\l-\mu)(\l-\nu)}\left(\Xi -(\mu+\nu)H+\mu\nu E\right)
\\
q^2&=\frac1{(\mu-\nu)(\mu-\l)}\left(\Xi -(\nu+\l)H+\nu\l E\right)
\\
r^2&=\frac1{(\nu-\l)(\nu-\mu)}\left(\Xi -(\l+\mu)H+\l\mu E\right)
\end{aligned}
\end{equation}
so that
$$
\begin{aligned}
(\mu-\nu)q^2r^2&=-\frac{\left(\Xi -(\nu+\l)H+\nu\l E\right)\left(\Xi -(\l+\mu)H+\l\mu E\right)}{(\l-\mu)(\l-\nu)(\mu-\nu)}
\\
(\nu-\l)r^2p^2&=-\frac{\left(\Xi -(\l+\mu)H+\l\mu E\right)\left(\Xi -(\mu+\nu)H+\mu\nu E\right)}{(\l-\mu)(\l-\nu)(\mu-\nu)}
\\
(\l-\mu)p^2q^2&=-\frac{\left(\Xi -(\mu+\nu)H+\mu\nu E\right)\left(\Xi -(\nu+\l)H+\nu\l E\right)}{(\l-\mu)(\l-\nu)(\mu-\nu)}
\end{aligned}
$$
Later on, we shall use the notations
\begin{equation}
\label{DfntRoots}
\begin{aligned}
x_-(\l,\mu,\nu)&=(\mu+\nu)H-\mu\nu E
\\
x_0\,(\l,\mu,\nu)&=(\mu+\l)H-\mu\l E
\\
x_+(\l,\mu,\nu)&=(\l+\nu)H-\l\nu E
\end{aligned}
\end{equation}
Therefore, we find that $\Xi $ satisfies the second order ODE
$$
\begin{aligned}
\ddot{\Xi }=&-2K_{\l,\mu,\nu}\left((\Xi -x_-(\l,\mu,\nu))(\Xi -x_0(\l,\mu,\nu))\right.
\\
&\left.+(\Xi -x_0(\l,\mu,\nu))(\Xi -x_+(\l,\mu,\nu))+(\Xi -x_+(\l,\mu,\nu))(\Xi -x_0(\l,\mu,\nu))\right)
\end{aligned}
$$
which can be put in the form
\begin{equation}
\label{ODEOm}
\ddot\Xi =-2K_{\l,\mu,\nu}P'_{\l,\mu,\nu}(\Xi)
\end{equation}
where $P_{\l,\mu,\nu}$ is the cubic
\begin{equation}
\label{DfntP}
P_{\l,\mu,\nu}(X)=(X-x_-(\l,\mu,\nu))(X-x_0(\l,\mu,\nu))(X-x_+(\l,\mu,\nu))
\end{equation}
and
\begin{equation}
\label{DfntK}
K_{\l,\mu,\nu}=\frac{\l^2(\mu-\nu)+\mu^2(\nu-\l)+\nu^2(\l-\mu)}{(\l-\mu)(\l-\nu)(\mu-\nu)}
\end{equation}
In the sequel, we assume that the initial data for $(p,q,r)$ is such that
$$
r(0)=0\,,\qquad p(0)(q(0)\not=0
$$
Let us compute
\begin{equation}
\label{FrmlRoots}
\begin{aligned}
x_-(\l,\mu,\nu)&=\l\nu p(0)^2+\mu^2q(0)^2+\mu(\l-\nu)p(0)^2
\\
x_0\,(\l,\mu,\nu)&=\l^2 p(0)^2+\mu^2q(0)^2
\\
x_+(\l,\mu,\nu)&=\l^2 p(0)^2+\left(\frac{\nu+\l}{\mu}-\frac{\nu\l}{\mu^2}\right)\mu^2q(0)^2
\end{aligned}
\end{equation}
We shall also assume that
\begin{equation}
\label{AssEigV}
\l>\mu>\nu\,,\qquad\l\nu<0\,,\qquad |\mu|\ll \l\hbox{ and }\l\sim |\nu|
\end{equation}
Then $K_{\l,\mu,\nu}>0$ --- in fact $K_{\l,\mu,\nu}\sim 2$, and $\Xi$ is a periodic function
of $t$ such that
\begin{equation}
\label{MinMaxOm}
\inf_{t\in\bR}\Xi (t)=x_0(\l,\mu,\nu)\,,\quad\sup_{t\in\bR}\Xi (t)=x_+(\l,\mu,\nu)
\end{equation}
with half-period
\begin{equation}
\label{DfntT}
T_{\l,\mu,\nu}=\frac1{2\sqrt{K_{\l,\mu,\nu}}}
    \int_{x_0(\l,\mu,\nu)}^{x_+(\l,\mu,\nu)}\frac{dx}{\sqrt{-P_{\l,\mu,\nu}(x)}}
\end{equation}
We are interested in the growth of the (squared) $\mathbf{H}^3$ norm
\begin{equation}
\label{DfntW}
W(t)=\l^6p(t)^2+\mu^6q(t)^2+\nu^6r(t)^2
\end{equation}
Expressing $p^2$, $q^2$ and $r^2$ in terms of $E$, $H$ and $\Xi $, it is found that
\begin{equation}
\label{EHOM>W}
W=\frac{\l^6(\Xi-x_-(\l,\mu,\nu))}{(\l-\mu)(\l-\nu)}
+\frac{\mu^6(\Xi-x_+(\l,\mu,\nu))}{(\mu-\nu)(\mu-\l)}
+\frac{\nu^6(\Xi-x_0(\l,\mu,\nu))}{(\nu-\l)(\nu-\mu)}
\end{equation}
Hence, when $\Xi=x_+(\l,\mu,\nu)$, then
$$
\begin{aligned}
W&=\frac{\l^6(x_+(\l,\mu\,\nu)-x_-(\l,\mu,\nu))}{(\l-\mu)(\l-\nu)}
+\frac{\nu^6(x_+(\l,\mu\,\nu)-x_0(\l,\mu,\nu))}{(\nu-\l)(\nu-\mu)}
\\
&\ge\frac{\l^6(x_+(\l,\mu\,\nu)-x_-(\l,\mu,\nu))}{(\l-\mu)(\l-\nu)}
\end{aligned}
$$
Let us compute
\begin{equation}
\label{x+-x-}
\begin{aligned}
x_+(\l,\mu\,\nu)-x_-(\l,\mu,\nu)&=
(\l-\mu)(\l-\nu)p(0)^2+\left(\frac{\nu+\l}{\mu}-\frac{\nu\l}{\mu}-1\right)\mu^2q(0)^2
\\
&\gtrsim-\nu\l q(0)^2\sim\l^2q(0)^2
\end{aligned}
\end{equation}
We shall pick the initial data such that
\begin{equation}
\label{W(0)}
W(0)=\l^6p(0)^6+\mu^6q(0)^6\hbox{ with }\l^6p(0)^2\sim\tfrac12 W(0)\hbox{ and }\mu^6q(0)^2\sim\tfrac12W(0)
\end{equation}
Hence, when $\Xi$ reaches $x_+(\l,\mu,\nu)$, one has
\begin{equation}
\label{Wmax}
W\gtrsim\frac{\l^8q(0)^2}{(\l-\mu)(\l-\nu)}\sim \tfrac12\frac{\l^8}{\mu^6(\l-\mu)(\l-\nu)}W(0)
\sim\tfrac14\frac{\l^6}{\mu^6}W(0)\,.
\end{equation}
Hence $W$ jumps from $W(0)$ to a quantity $\sim\tfrac14\frac{\l^6}{\mu^6}W(0)$ in an interval of time that
does not exceed one period of the $\Xi$ motion, i.e. $2T_{\l,\mu,\nu}$. Let us estimate this interval of
time. We recall the asymptotic equivalent for the period of an elliptic integral in the modulus 1 limit.

\begin{lem}
Assume that $x_-<x_0<x_+$. Then
$$
\int_{x_0}^{x_+}\frac{dx}{\sqrt{(x-x_-)(x-x_0)(x_+-x)}}
\sim\frac1{\sqrt{x_+-x_-}}\ln\left(\frac1{1-\sqrt{\frac{x_+-x_0}{x_+-x_-}}}\right)
$$
uniformly in $x_-$, $x_0$, and $x_+$ as $\frac{x_+-x_0}{x_+-x_-}\to 1$.
\end{lem}

Here
$$
\frac1{\sqrt{x_+(\l,\mu,\nu)-x_-(\l,\mu,\nu)}}\lesssim\frac1{\sqrt{\l^2q(0)^2}}
    \sim\frac{|\mu|^3}{\l}\sqrt{\frac2{W(0)}}
$$
Next
\begin{equation}
\label{x0-x-}
x_0(\l,\mu,\nu)-x_-(\l,\mu,\nu)=(\l-\mu)(\l-\nu)p(0)^2
\end{equation}
so that
$$
\begin{aligned}
\frac1{1-\sqrt{\frac{x_+-x_0}{x_+-x_-}}}&\sim
\frac1{1-\sqrt{1-\frac{(\l-\mu)(\l-\nu)p(0)^2}{(\l-\mu)(\l-\nu)p(0)^2+(\mu(\nu+\l)-\nu\l-\mu^2)q(0)^2}}}
\\
&\sim\frac{(\l-\mu)(\l-\nu)p(0)^2+(\mu(\nu+\l)-\nu\l-\mu^2)q(0)^2}{2(\l-\mu)(\l-\nu)p(0)^2}
\\
&\sim\frac{q(0)^2}{2p(0)^2}\sim\tfrac12\frac{W(0)/2\mu^6}{W(0)/2\l^6}=\frac{\l^6}{2\mu^6}
\end{aligned}
$$
Hence
\begin{equation}
\label{Tsim}
2T_{\l,\mu,\nu}\lesssim\frac2{\sqrt{W(0)}}\frac{\mu^3}{\l}\ln\left(\frac{\l^6}{2\mu^6}\right)
    \le \frac{12}{\sqrt{W(0)}}\frac{|\mu|^3}{\l}\ln\left(\frac{\l}{\mu}\right)
\end{equation}
{\bf Conclusion:} collecting (\ref{W(0)}), (\ref{Wmax}) and (\ref{Tsim}), we see that the squared
$\mathbf{H}^3$ norm $W$ varies from $W(0)$ to a quantity $\sim\rho^6W(0)$ in an interval of time $\lesssim
\frac{12}{\sqrt{W(0)}}\frac{\mu^2\ln\rho}{\rho}$. (Here $\rho=\l/\mu$).

%New typing starting Jan. 22, 2007---------

We now proceed to obtain similar bursting estimates for the enstrophy. We return to
 (\ref{MinMaxOm}) and (\ref{DfntT}). Pick the initial data so that
\[
    \Xi (0)=\lambda^2p(0)^2+\mu^2q(0)^2 \mbox{ with } 
   \lambda^2p(0)^2 \sim \frac{1}{2}\Xi (0) \mbox{ and } 
   \mu^2q(0)^2\sim \frac{1}{2}\Xi (0). 
\]
Then
$$
\begin{aligned}
   x_+(\lambda , \mu , \nu )&-x_{-}(\lambda , \mu , \nu )
  \\
  &=(\lambda -\mu ) (\lambda -\nu ) p(0)^2 +
               \left( \frac{\nu +\lambda }{\mu}-\frac{\nu\lambda}{\mu^2}-1\right) \mu^2 q(0)^2 
  \\
   &\sim 2 \lambda^2 p(0)^2+\lambda^2 q(0)^2 \sim \frac{\lambda^2}{\mu^2}\Xi (0)
\end{aligned}
$$
while
\[ x_0(\lambda , \mu , \nu )-x_{-}(\lambda , \mu , \nu ) =
   (\lambda -\mu )(\lambda -\nu ) p(0)^2\sim 2\lambda^2 p(0)^2 \sim \Xi (0). 
\]
Hence, in the limit as $\rho = \lambda / |\mu| \rightarrow+\infty$, one has
$$
\begin{aligned} 
 2T_{\lambda , \mu , \nu}&\sim\frac{1}{2\sqrt{2}}\frac{1}{\sqrt{\rho^2\Xi (0)}} 
    \ln \frac{1}{1-\sqrt{1-\frac{\Xi (0)}{\frac{1}{2}\rho^2\Xi (0)} } } 
\\
&= \frac{1}{2\sqrt{2\Xi (0)}}\frac{1}{\rho} 
   \ln \frac{1}{1-\sqrt{1-2\rho^{-2}}}
\sim\frac{1}{\sqrt{2\Xi (0)}} \frac{\ln \rho}{\rho}. 
\end{aligned}
$$ 
And $\Xi$ varies from
\[ x_0(\lambda , \mu , \nu )=\Xi (0) \mbox{ to } x_{+}(\lambda , \mu , \nu )
   \sim \rho^2 \Xi (0)
\]
on an interval of time  of length $T_{\lambda , \mu , \nu}$. $\ \ \blacksquare$

%--------------------------------------------------------
  \section{Strictly resonant Euler systems: the case of 3-waves resonances on
      small-scales}
%---------------------------------------------------------

%--------------------------------------------
 \subsection{ Infinite dimensional uncoupled $SO(3)$ systems }
%---------------------------------------------

In this section, we consider the 3-wave resonant set $K^*$ when 
\[ |k|^2,\ |m|^2,\ |n|^2\geq \frac{1}{\eta^2},\ 0<\eta\ll 1, \]
 i.e. 3-wave resonances 
on small scales; here $|k|^2=k_1^2+k_2^2+k_3^2$, where
$(k_1, k_2, k_3)$ index the curl eigenvalues, and
similarly for $|m|^2,\ |n|^2$. Recall that $k_2+m_2=n_2,\
k_3+m_3=n_3$ (exact convolutions), but  that the summation on $k_1,\ m_1$ on the right hand
side of Eqs. (\ref{resonant_euler_again}) is not a convolution. However, for $|k|^2,\ |m|^2,\ |n|^2\geq \frac{1}{\eta^2},$ the summation in $k_1,\ m_1$ becomes 
an asymptotic convolution. First:
\begin{prop}\label{prop4.1} 
  The set $K^*$ restricted to $|k|^2, |m|^2, |n|^2\geq \frac{1}{\eta^2}$,\ 
  $\forall \eta, 0<\eta\ll 1$ is not empty: there exist at least
one $h/R$ with resonant three waves satisfying the above small scales condition.
\end{prop}

\textit{Proof.} We follow the algebra of the exact transcendental dispersion law
(\ref{trancendental}) derived in the proof of Lemma \ref{lemma5.7}.
Note that $\tilde{P}(\vartheta_3)<0$ for $\vartheta_3=\frac{1}{h^2}$
large enough. We can choose $h_m=\frac{\beta^2(m_1,m_2,\alpha m_3)}{m_3^2}=0,$
say in the specific limit
$\frac{h}{2\pi m_3}\rightarrow 0$, and $\beta (m_1, m_2, \alpha m_3)\sim m_1\pi 
+m_2\frac{\pi}{2} +\frac{\pi}{4}$. Then $\tilde{P}_0=h_k^2h_n^2 >0$ and
$\tilde{P}(\vartheta_3)$ must possess at least one (transcendental) root
$\vartheta_3=\frac{1}{h^2}.\ \ \ \blacksquare$
\medskip

In the above context, the radial components of the curl eigenfunctions  involve
cosines and sines in $\frac{\beta r}{R}$ (cf. Section 3, \cite{M-N-B-G})
and the summation in $k_1,\ m_1$ on the right hand side of the resonant Euler equations 
(\ref{resonant_euler_again}) becomes an \textit{asymptotic convolution}.
The rigorous asymptotic convolution estimates are highly technical and
detailed in \cite{Fro-M-N}. The 3-wave resonant systems
for $|k|^2,\ |m|^2,\ |n|^2\geq \frac{1}{\eta^2}$ are equivalent to those
of an equivalent periodic lattice
$[0, 2\pi ]\times [0, 2\pi ]\times [0, 2\pi h ],\ \vartheta_3=\frac{1}{h^2}$;
the resonant three wave relation becomes:

\begin{subequations}\label{triplets}
  \begin{align}
 \pm \left( \vartheta_3+\vartheta_1\frac{n_1^2}{n_3^2}+
     \vartheta_2\frac{n_2^2}{n_3^2}\right)^{-\frac{1}{2}}
  &\pm \left( \vartheta_3+\vartheta_1\frac{k_1^2}{k_3^2}+
     \vartheta_2\frac{k_2^2}{k_3^2}\right)^{-\frac{1}{2}} \nonumber\\
  &\pm \left( \vartheta_3+\vartheta_1\frac{m_1^2}{m_3^2}+
     \vartheta_2\frac{m_2^2}{m_3^2}\right)^{-\frac{1}{2}}  = 0,\\
  &k+m=n,\ k_3m_3n_3 \neq 0. \hspace{1in}
  \end{align}
\end{subequations}
The algebraic geometry of these rational 3-wave resonance equations has
been investigated in depth in \cite{B-M-N3} and
\cite{B-M-N4}. Here $\vartheta_1,\ \vartheta_2,\ \vartheta_3$ are periodic lattice
parameters;
in the small-scales cylindrical case, $\vartheta_1=\vartheta_2=1$ (after rescaling
of $n_2$, $k_2$, $m_2$), \ $\vartheta_3=1/h^2,\ 
h$ height.
Based on the algebraic geometry of ``resonance curves'' in \cite{B-M-N3},
\cite{B-M-N4},
we investigate the resonant 3D Euler equations (\ref{resonant_euler_again})
in the equivalent \textit{periodic lattices}.

First, triplets $(k,m,n)$ solution of (\ref{triplets}) are invariant 
under the reflection symmetries $\sigma_0, \sigma_1, \sigma_2, \sigma_3$ defined 
in Corollary \ref{cor3.1} and Remark \ref{rem3.2}:
$\sigma_0=Id, \, \sigma_j(k)=(\epsilon_{i,j}k_i),\ 1\leq i\leq 3,\
\epsilon_{i,j}=+1$ if $ i\neq j$, $\epsilon_{i,j}=-1$ if $i=j$, $1\leq j\leq 3$.
Second the set $K^*$ in (\ref{triplets}) is invariant under the homothetic transformations:
\begin{equation}
  (k, m, n) \rightarrow (\gamma k,\gamma m,\gamma n),\ \gamma \mbox{ rational}.
\end{equation}
The resonant triplets lie on projective lines in the wavenumber space, with
equivariance under $\sigma_j,\ 0\leq j\leq 3$ and $\gamma$-rescaling. For every given equivariant family of such
projective lines, the resonant curve is the graph of $\frac{\vartheta_3}{\vartheta_1}$ versus
$\frac{\vartheta_2}{\vartheta_1}$, for parametric domain resonances in
$\vartheta_1,\ \vartheta_2,\ \vartheta_3$.

\begin{lem}\label{lem4.2} (p.17, {\em \cite{B-M-N4}}). For every equivariant $(k,m,n)$,
the resonant curve in the quadrant $\vartheta_1 >0, \vartheta_2>0, \vartheta_3>0$ is the graph of
a smooth function $\vartheta_3/\vartheta_1\equiv F(\vartheta_2/\vartheta_1)$ intersected
with the quadrant.
  
\end{lem}

\begin{thm}\label{thm4.2}(p.19, {\em \cite{B-M-N4}}). A resonant curve in the quadrant
{\small $\vartheta_3/\vartheta_1>0,$ $\vartheta_2/\vartheta_1>0$} is called irreducible if:
\begin{equation}\label{irr_det}
  det \left( \begin{array}{ccc}
                  k_3^2 & k_2^2 & k_1^2 \\
                  m_3^2 & m_2^2 & m_1^2 \\
                  n_3^2 & n_2^2 & n_1^2 
             \end{array}
       \right) \neq 0.
\end{equation}
An irreducible resonant curve is uniquely characterized by six non-negative algebraic
invariants $\mathcal{P}_1,\ \mathcal{P}_2,\ \mathcal{R}_1,
\ \mathcal{R}_2,\ \mathcal{S}_1,\ \mathcal{S}_2$, such that
\begin{displaymath}
\left\{ \frac{n_1^2}{n_3^3}, \frac{n_2^2}{n_3^2}\right\} = 
 \left\{ \mathcal{P}_1^2, \mathcal{P}_2^2\right\} ,\ \ 
\left\{ \frac{k_1^2}{k_3^3}, \frac{k_2^2}{k_3^2}\right\} = 
 \left\{ \mathcal{R}_1^2, \mathcal{R}_2^2\right\} ,\ \ 
\left\{ \frac{m_1^2}{m_3^3}, \frac{m_2^2}{m_3^2}\right\} = 
 \left\{ \mathcal{S}_1^2, \mathcal{S}_2^2\right\} ,\ \ 
\end{displaymath}
and permutations thereof.
\end{thm}

\begin{lem}\label{lem4.3} (p. 25, {\em \cite{B-M-N4}}). For resonant triplets $(k,m,n)$
associated to a given irreducible resonant curve, that is verifying Eq. 
(\ref{irr_det}), consider the convolution equation $n=k+m$. 
Let $\sigma_i(n)\neq n,\ \forall i,\ 1\leq i\leq 3$. Then there are no more that \textit{two}
solutions $(k,m)$ and $(m,k)$, for a given $n$, provided the six 
 non-degeneracy conditions (3.39)-(3.44) in {\em \cite{B-M-N4}} 
for the algebraic invariants of the irreducible curve are verified.
  
\end{lem}

For more details on the technical non-degeneracy conditions,
see the Appendix. An exhaustive algebraic geometric investigation of all 
solutions to $n=k+m$ on irreducible resonant curves is found in \cite{B-M-N4}.
The essence of the above lemma lies in that given such an irreducible,
``non-degenerate'' triplet $(k,m,n)$ on $K^*$, all other triplets
on the same irreducible resonant curves are exhaustively given by the equivariant 
projective lines:
\begin{eqnarray}
     (k,m,n) &\rightarrow& (\gamma k, \gamma m, \gamma n),\ \ \mbox{ for some } \gamma
              \mbox{ rational }, \\
     (k,m,n) &\rightarrow& (\sigma_j k, \sigma_j m, \sigma_j n), \ \ j=1,2,3, 
\end{eqnarray}
and permutations of $k$ and $m$ in the above.
Of course the homothety $\gamma$ and the $\sigma_j$ symmetries preserve
the convolution. This context of irreducible, ``non-degenerate'' resonant curves
yields an infinite dimensional, uncoupled system of rigid body $SO(3;\mathbf{R})$ and 
$SO(3;\mathbf{C})$ dynamics for the 3D resonant Euler equations  (\ref{resonant_euler_again}).
\begin{thm}\label{thm4.4}
   For any irreducible triplet $(k,m,n)$ which satisfy Theorem \ref{thm4.2}, and
   under the ``non-degeneracy'' conditions of Lemma \ref{lem4.3} (cf. Appendix),
the resonant Euler equations split into the infinite, countable sequence of 
uncoupled $SO(3;\mathbf{R})$ systems:
\begin{subequations}\label{eq4.7}
   \begin{align}
      \dot{a}_k = \Gamma_{kmn}(\lambda_m-\lambda_n) a_ma_n, \\
      \dot{a}_m = \Gamma_{kmn}(\lambda_n-\lambda_k) a_na_k,\hspace{.05in} \\
      \dot{a}_n = \Gamma_{kmn}(\lambda_k-\lambda_m) a_ka_m,
   \end{align}
\end{subequations}
\begin{equation}\label{eq4.8.1}
 \mbox{for all } (k,m,n)=\gamma (\sigma_j(k^*), \sigma_j(m^*),\sigma_j(n^*)),\ 
\gamma = \pm 1, \pm 2, \pm 3...,  \ 0\leq j\leq 3.
\end{equation}
$k^*, m^*, n^*$ are some relatively prime integer vectors in $\mathbf{Z}^3$
characterizing the equivariant family of projective lines $(k,m,n)$;
$\Gamma_{kmn}=i <\mathbf{\Phi}_k\times\mathbf{\Phi}_m, \mathbf{\Phi}_n^*>,$ 
$\Gamma_{kmn}$ real.
\end{thm}

\textit{Proof.} Theorem \ref{thm4.4}  is a simpler version for invariant manifolds of
more general $SO(3;\mathbf{C})$ systems. It is a straightforward corollary of Proposition 3.2,
Proposition 3.3, Theorem 3.3, Theorem 3.4 and Theorem 3.5 in \cite{B-M-N4}.
The latter article did not explicit the resonant equations and did not use 
the curl-helicity algebra fundamentally underlying this present work.
Rigorously asymptotic infinite countable sequences of uncoupled 
$SO(3;\mathbf{R}),\ SO(3;\mathbf{C})$ systems are not derived via
the usual harmonic analysis tools of Fourier modes, in the 3D Euler context.
Polarization of curl eigenvalues and eigenfunctions and helicity play an essential role.

\begin{cor}
  Under the conditions $\lambda_{n^*}-\lambda_{k^*}>0,\ \lambda_{k^*}-\lambda_{m^*}>0,$
 the resonant Euler systems (\ref{eq4.7}) admit a disjoint, countable family of homoclinic
  cycles. Moreover, under the conditions $\lambda_{n^*}\gg +1,\ \lambda_{m^*}\ll -1,\ 
  |\lambda_{k^*}|\ll \lambda_{n^*},$ 
 each subsystem (\ref{eq4.7}) possesses orbits
whose $\mathbf{H}^s$ norms, $s\geq 1$, burst
arbitrarily large in arbitrarily small times.
\end{cor}

\begin{rem}
  One can prove that there exists some $\Gamma_{max},\ 0<\Gamma_{max}<\infty ,$
such that $|\Gamma_{kmn}|<\Gamma_{max}$, for all $(k,m,n)$ on the equivariant 
projective lines defined by (\ref{eq4.8.1}). Systems (\ref{eq4.7}) ``freeze''
cascades of energy; their total enstrophy
$\Xi (t)=\sum_{(k,m,n)} (\lambda_k^2a_k^2(t)+\lambda_m^2a_m^2(t)+\lambda_n^2a_n^2(t) )$
remains bounded, albeit with large bursts of
$\Xi (t)/\Xi (0)$, on the reversible orbits  topologically close to the homoclinic cycles.
\end{rem}

%---------------------------------------------------
\subsection{Coupled $SO(3)$ rigid body resonant systems}
%----------------------------------------------------

We now derive a new resonant Euler system which couples
{\bf two} $SO(3;\mathbf{R})$ rigid bodies via a common principle axis of
inertia and a common moment of inertia. This 5-dimensional system
conserves energy, helicity, and is rather interesting in that dynamics on its
homoclinic manifolds show bursting cascades of enstrophy to the smallest
scale in the resonant set.
We consider the equivalent periodic lattice geometry under the conditions of
Proposition \ref{prop4.1}.

In Appendix, we prove that for an ``irreducible'' 3-wave resonant set
which now satisfies the algebraic ``degeneracy'' (\ref{dege}), there exist exactly 
two ``primitive'' resonant triplets $(k,m,n)$ and $(\tilde{k},\tilde{m},n)$, 
where $k,\ m,\ \tilde{k},\ \tilde{m}$ are relative prime integer valued vectors
in $\mathbf{Z}^3$:
\begin{lem}\label{lem4.8}
  Under the algebraic degeneracy condition (\ref{dege}) the irreducible equivariant
family of projective lines
in $K^*$ is exactly generated by the following two ``primitive'' triplets:
  \begin{subequations}\label{eq4.8}
     \begin{align}
        n=k+m,\ \ k=a\overline{k},\ \ m=b\overline{m}, \hspace{.8in}\\
        n=\tilde{k}+\tilde{m}, \ \ 
        \tilde{k}=a^\prime \sigma_i(\overline{k})+b^\prime\sigma_j(\overline{m}),\hspace{.55in} \\
        \mbox{that is,} \hspace{3.6in}\nonumber \\
        n=a\overline{k}+b\overline{m}, \hspace{1.85in}\\
        n=a^\prime \sigma_i(\overline{k})+b^\prime\sigma_j(\overline{m}),\hspace{1.3in}
     \end{align}
  \end{subequations}
where $\sigma_i\neq\sigma_j$ are some reflection symmetries,
$a,\ b,\ a^\prime,\ b^\prime$ are relatively prime integers, positive or negative, and
$\overline{k},\ \overline{m}$ are relatively prime integer valued vectors in
$\mathbf{Z}^3$, that is:
\[ (a, a^\prime )=(b, b^\prime)=(a,b)=(a^\prime ,b^\prime )=1, \ \ 
   (\overline{k},\overline{m})=1,
\]
where $(\ ,\ )$ denotes the Greatest Common Denominator of two integers. 
All other resonant wave number triplets are generated by the group actions
$\sigma_l,$ \ $l=1,2,3,$ and homothetic rescalings 
$(k,m,n)\rightarrow \gamma (k,m,n), $\
 $(\tilde{k},\tilde{m},n)\rightarrow \gamma (\tilde{k},\tilde{m},n),$
$ (\gamma\in \mathbf{Z})$
of the ``primitive'' triplets.
\end{lem}
\begin{rem}
  It can be proven that the set of such coupled ``primitive'' triplets is 
not empty on the periodic lattice. The algebraic  irreducibility condition
of Lemma \ref{lem4.2} implies that $\pm k_3/|k|=\pm \tilde{k}_3/|\tilde{k}|$
and $\pm m_3/|m|=\pm \tilde{m}_3/|\tilde{m}|$, which is obviously
verified in equations (\ref{eq4.8}).
\end{rem}
\begin{thm}\label{Theorem4.10}
  Under conditions of Lemma \ref{lem4.8} the resonant Euler system reduces to
a system of two rigid bodies coupled via $a_n(t)$:
\begin{subequations}\label{resonant_system}
   \begin{align}
     \dot{a}_k&=(\lambda_m-\lambda_n)  \Gamma a_m a_n \\
      \dot{a}_m&=(\lambda_n-\lambda_k)  \Gamma a_n a_k \\
     \dot{a}_n&=(\lambda_k-\lambda_m)  \Gamma a_k a_m 
            + (\lambda_{\tilde{k}}-\lambda_{\tilde{m}})\tilde{\Gamma}a_{\tilde{k}}a_{\tilde{m}}\\
     \dot{a}_{\tilde{m}}&=(\lambda_n-\lambda_{\tilde{k}})\tilde{\Gamma}a_n a_{\tilde{k}} \\
     \dot{a}_{\tilde{k}}&=(\lambda_{\tilde{m}}-\lambda_n)\tilde{\Gamma}a_{\tilde{m}} a_n, 
   \end{align}
\end{subequations}
where $\Gamma = i <\mathbf{\Phi}_k\times \mathbf{\Phi}_m, \mathbf{\Phi}_n^*>,\ 
  \tilde{\Gamma}=i <\mathbf{\Phi}_{\tilde{k}}\times\mathbf{\Phi}_{\tilde{m}}, \mathbf{\Phi}_n^*>$.
Energy and Helicity are conserved.
\end{thm}

\begin{thm}
  The resonant system (\ref{resonant_system}) possesses three independent
conservation laws:
\begin{subequations}\label{resonant_system/other}
  \begin{align}
    \mathcal{E}_1 &= a_k^2+(1-\alpha )a_m^2, \label{eq4.10a}\\
    \mathcal{E}_2 &=a_n^2+\alpha a_m^2 + (1-\tilde{\alpha})a_{\tilde{m}}^2,\label{eq4.10b} \\
    \mathcal{E}_3 &=a_{\tilde{k}}^2+\tilde{\alpha}a_{\tilde{m}}^2, \label{eq4.10c}
  \end{align}
\end{subequations}
where
\begin{subequations}
   \begin{align}
      \alpha = (\lambda_m-\lambda_k)/(\lambda_n-\lambda_k), \\
      \tilde{\alpha} = (\lambda_{\tilde{m}}-\lambda_n)/(\lambda_{\tilde{k}}-\lambda_n).
   \end{align}
\end{subequations}
\end{thm}
\begin{thm}
  Under the conditions
  \begin{subequations}\label{eq4.12}
    \begin{align}
      \lambda_m < \lambda_k < \lambda_n, \\
       \lambda_{\tilde{m}} < \lambda_n < \lambda_{\tilde{k}},
    \end{align}
  \end{subequations}
which imply $\alpha <0,\ \tilde{\alpha}<0,$ the equilibria
$(\pm a_k(0), 0, 0, 0, \pm a_{\tilde{k}}(0))$ are hyperbolic
for $|a_{\tilde{k}}(0)|$ small enough with respect to
$|a_k(0)|$. The unstable manifolds of these equilibria are 
one dimensional, and the nonlinear dynamics of system (\ref{resonant_system})
are constrained on the ellipse $\mathcal{E}_1$ (\ref{eq4.10a})
for $a_k(t)$, $a_m(t)$, the hyperbola $\mathcal{E}_3$ (\ref{eq4.10c}) for
$a_{\tilde{k}}(t)$, $a_{\tilde{m}}(t)$, and the hyperboloid $\mathcal{E}_2$ 
(\ref{eq4.10b}) for
$a_m(t)$, $a_{\tilde{m}}(t)$, $a_n(t)$.
\end{thm}
\begin{thm}\label{thm4.13}
Let the 2-manifold $\mathcal{E}_1\cap\mathcal{E}_2\cap\mathcal{E}_3$ be
coordinatized by  $(a_m, a_{\tilde{m}})$. On this 2-manifold, the resonant 
system (\ref{resonant_system}) is Hamiltonian, and therefore integrable.
Its Hamiltonian vector field $\mathbf{h}$ is defined by 
\begin{equation}
\iota_\mathbf{h}\omega\ =\ \Gamma(\lambda_n-\lambda_k)\frac{da_{\tilde m}}{a_{\tilde k}}
-
\tilde\Gamma(\lambda_n-\lambda_{\tilde k})\frac{da_m}{a_k},
\end{equation}
where $\iota_\mathbf{h}\omega$ designates the inner product  of the 
symplectic 2-form 
\begin{equation}
\omega=\frac{da_m\wedge da_{\tilde m}}{a_ka_na_{\tilde k}}
\end{equation}
with the vector field $\mathbf{h}$.
\end{thm}

\textit{Proof of Theorem \ref{thm4.13}}:  Eliminating $a_k(t)$ via $\mathcal{E}_1$,
  $a_n(t)$ via $\mathcal{E}_2$, $a_{\tilde{k}}(t)$ via $\mathcal{E}_3$,
the resonant system (\ref{resonant_system}) reduces to:
\begin{eqnarray*}
   \dot{a}_m &=& \pm \Gamma (\lambda_n-\lambda_k) 
                   (\mathcal{E}_1-(1-\alpha )a_m^2)^{\frac{1}{2}}
                   (\mathcal{E}_2-\alpha a_m^2+(\tilde{\alpha}-1)a_{\tilde{m}}^2)^{\frac{1}{2}}\\
   \dot{a}_{\tilde{m}} &=& \pm \tilde{\Gamma} (\lambda_n-\lambda_{\tilde{k}}) 
                   (\mathcal{E}_2-\alpha a_m^2+(\tilde{\alpha}-1)a_{\tilde{m}}^2)^{\frac{1}{2}}
                   (\mathcal{E}_3-\tilde{\alpha}
                         a_{\tilde{m}}^2)^\frac{1}{2};
\end{eqnarray*}
after changing the time variable into
\[ t\rightarrow\int_0^t(\mathcal{E}_1-(1-\alpha )a_m^2)^{\frac{1}{2}}
    (\mathcal{E}_2-\alpha a_m^2+(\tilde{\alpha}-1)a_{\tilde{m}}^2)^{\frac{1}{2}}
    (\mathcal{E}_3-\tilde{\alpha}a_{\tilde{m}}^2)^{\frac{1}{2}}ds\,.
\]

On each component of the manifold $\mathcal{E}_1\cap\mathcal{E}_2\cap\mathcal{E}_3$,
the following functionals are conserved:

\begin{eqnarray*}
     \mathcal{H}(a_m, a_{\tilde{m}}) = &
      \pm &\tilde{\Gamma} (\lambda_n-\lambda_{\tilde k})
        \int \frac{da_{m}}{(\mathcal{E}_1-(1-\alpha )a_m^2)^{1/2}} \hspace{.7in}\\
      &\pm&\Gamma (\lambda_n-\lambda_{k})
        \int \frac{da_{\tilde{m}}}{(\mathcal{E}_3-\tilde{\alpha} a_{\tilde{m}}^2)^{1/2}}. 
        \hspace{.7in} \ \ \ \ \ \blacksquare
\end{eqnarray*}

\medskip

\noindent 
Observe that the system of two coupled rigid bodies (\ref{resonant_system}) 
does not seem to admit a simple Lie-Poisson bracket in the original variables 
$(a_k,a_m,a_n,a_{\tilde m},a_{\tilde k})$. Yet, when 
restricted to the 2-manifold $\mathcal{E}_1\cap\mathcal{E}_2\cap\mathcal{E}_3$
that is invariant under the flow of (\ref{resonant_system}), it is Hamiltonian 
and therefore integrable. 

This raises the following interesting issue: according to the shadowing 
Theorem \ref{theorem2.10}, the Euler dynamics remains asymptotically close 
to that of chains of coupled $SO(3;\mathbf{R})$ and $SO(3;\mathbf{C})$ rigid 
body systems. Perhaps some new information could be obtained in this way.
We are currently investigating this question and will report on it in a forthcoming
publication \cite{GMNNew}.

Already the simple 5-dimensional system (\ref{resonant_system}) has interesting
dynamical properties, wich we could not find in the existing literature on systems 
related to spinning tops.

Consider for instance the dynamics of the resonant system (\ref{resonant_system}) 
with I.C. topologically close to the hyperbola equilibria
$(\pm a_k(0), 0,0,0, \pm a_{\tilde{k}}(0))$. 
Under the conditions of (\ref{eq4.12}) and with the help of the integrability 
Theorem \ref{thm4.13}, it is easy to construct equivariant families
of homoclinic cycles at these hyperbolic critical points:

\begin{cor}
  The hyperbolic critical points $(\pm a_k(0), 0,0,0,\pm a_{\tilde{k}}(0))$ 
possess 1-dimensional homoclinic cycles on the cones
\begin{equation}\label{cones}
   a_n^2+(1-\tilde{\alpha})a_{\tilde{m}}^2 = -\alpha a_m^2
\end{equation}
with $\alpha <0,\ \tilde{\alpha}<0$.
\end{cor}

Note that these are genuine homoclinic cycles, NOT sums of heteroclinic connections.
Initial conditions for the resonant system (\ref{resonant_system}) are now chosen 
in a small neighborhood of these hyperbolic critical points, the corresponding orbits are topologically close to these cycles.
With the ordering:
\begin{subequations}
  \begin{align}
    \lambda_m < \lambda_k < \lambda_n, \hspace{.2in}\\
    |\lambda_k| \ll |\lambda_m|,\ |\lambda_k|\ll \lambda_n, \\
   \lambda_{\tilde{m}} < \lambda_n < \lambda_{\tilde{k}}, \hspace{.2in}\\
   |\lambda_{\tilde{m}}| \ll \lambda_{\tilde{k}}, \hspace{.3in}\\
   \lambda_{\tilde{k}}\gg \lambda_n,\hspace{.3in}
  \end{align}
\end{subequations}
which can be realized with 
$|\frac{a^\prime}{a}|\gg 1$ and $|\frac{b^\prime}{b}|\ll 1$ in the resonant triplets
(\ref{eq4.8}), we can demonstrate bursting dynamics akin to Theorem \ref{thm3.9} and
\ref{thm3.11} for enstrophy and $\mathbf{H}^s$ norms, $s\geq 2$.
The interesting feature is the maximization of $|a_{\tilde{k}}(t)|$ near
the turning points of the homoclinic cycles on the cones  (\ref{cones}).
This corresponds to transfer of energy to the smallest scale 
$\tilde{k}$,
$\lambda_{\tilde{k}}.$

In a publication in preparation,  we investigate infinite systems
of the coupled rigid bodies equations 
(\ref{resonant_system}).

%--------------------------------------------------------

\renewcommand{\theequation}{A-\arabic{equation}}
 % redefine the command that creates the equation no.
\setcounter{equation}{0}  % reset counter 
\section*{\small APPENDIX}  % use *-form to suppress numbering

We focus on a resonant wave number triplet $(n,k,m)\in (\mathbf{Z}^{*})^3$
verifying
\begin{itemize}
   \item the convolution relation
      \begin{equation}\label{conv}
          n=k+m,
      \end{equation}
   \item the resonant 3-wave resonance relation
       \begin{equation}\label{3wave}
       \begin{aligned}
          \pm \frac{n_3}{\sqrt{\vartheta_1 n_1^2+\vartheta_2 n_2^2+\vartheta_3 n_3^2}}
          &\pm \frac{k_3}{\sqrt{\vartheta_1 k_1^2+\vartheta_2 k_2^2+\vartheta_3 k_3^2}}
           \\
           &\pm \frac{m_3}{\sqrt{\vartheta_1 m_1^2+\vartheta_2 m_2^2+\vartheta_3 m_3^2}}
            = 0, 
      \end{aligned}
      \end{equation}
   \item the condition of ``non-catalyticity''
        \begin{equation}\label{noncat}
            k_3m_3n_3 \neq 0,
        \end{equation}
   \item and the degeneracy condition of \cite{B-M-N4} (see p26)
       \begin{equation}\label{dege}
          G_{i,j}^{ir}(k,m) = k_i n_j m_l + k_l m_j n_i = 0,
       \end{equation}
       where $(i, j, l)$ is a permutation of $(1, 2, 3)$.
\end{itemize}
Then, we know (see lemma 3.5 (2) of \cite{B-M-N4}) that the system of equations 
(\ref{noncat})-(\ref{dege}) for the unknown $k$ and $m$, given the vector $n$,
admits exactly 4 solutions in $\mathbf{Z}^3\times\mathbf{Z}^3$:
\[ (k, m), \ \ (m,k),\ \ (\tilde{k}, \tilde{m}),\ \ (\tilde{m}, \tilde{k}).
\]
Here $k$ and $m$ are the two vectors of the original resonant triplet, whereas
\[ \tilde{k}=\alpha \sigma_i(k),\ \ \tilde{m}=\beta \sigma_j(m)
\]
where
\[ \alpha = \frac{m_ik_l-m_lk_i}{m_ik_l+m_lk_i}\notin \{ 0, \pm 1\}
   \mbox{ and } \beta = \frac{m_lk_j-m_jk_l}{m_lk_j+m_jk_l}\notin \{ 0, \pm 1\}
\]
and where the symmetries $\sigma_i$ and $\sigma_j$ are defined by
\[ \sigma_i : u=(u_l)_{l=1,2,3} \rightarrow 
    \left( (-1)^{\delta_{il}}u_l\right)_{l=1,2,3}.
\]

One verifies that 
\[ \sigma_i^2 = \sigma_j^2 = Id, \ \ \sigma_i\sigma_j = \sigma_j\sigma_i = -\sigma_l.
\]
That is, the group generated by $\sigma_i$ and $\sigma_j$ is the Klein group
$\mathbf{Z}/2\mathbf{Z}\times \mathbf{Z}/2\mathbf{Z}$.

Let us first write the irrational numbers $\alpha$ and $\beta$ under the irreducible 
representation
\[ \alpha = \frac{a^\prime}{a},\ \ \beta = \frac{b^\prime}{b},\ \  
   \mbox{ with } a, a^\prime , b, b^\prime \in \mathbf{Z}^* \mbox{ and }
   (a, a^\prime )=(b,b^\prime ) =1,
\]
where $(\ ,\ )$ denotes the Greatest Common Denominator of the integer pair.
From $\tilde{k}\in\mathbf{Z}^3$, it follows that $a|a^\prime k$; but since $(a, a^\prime )=1$, 
the Euclid's lemma yields that $a|k$. Similarly, $b|m$. Now set
\[ \overline{k}=\frac{1}{a}k\in \mathbf{Z}^3,\ \ \overline{m}=\frac{1}{b}m \in \mathbf{Z}^3. 
\]
Hence the integer vector $n$ admits the two decompositions
\[ n=a\overline{k}+b\overline{m}=
    a^\prime \sigma_i(\overline{k})+b^\prime \sigma_j(\overline{m}).
\]
Since the function
\[ z \longmapsto \frac{z_3}{\sqrt{\vartheta_1 z_1^2 + \vartheta_2 z_2^2 + \vartheta_3 z_3^2} }
\]
is homogeneous of degree $0$, we see that within the resonance condition (\ref{3wave})
 we can replace 
each vector $k, m$ and $n$ by any colinear vectors - either integer or not.
Suppose now that there exists some positive integer $d\neq 1$ such that $d|\overline{k}$;
then $d|n$, so that by setting 
\[ n_0=\frac{1}{d}n,\ \ k_0=\frac{1}{d}\overline{k},\ \ 
     m_0=\frac{1}{d}\overline{m} \]
we finally obtain
\[ n_0=ak_0+bk_0=a^\prime \sigma_i(k_0)+b^\prime \sigma_j(m_0).
\]
The triplets $(n_0, ak_0, bm_0)$ and $(n_0, a^\prime \sigma_i(k_0), b^\prime \sigma_j(m_0))$
further verify from the above remark, the convolution relation (\ref{conv}) and
the resonance relation (\ref{3wave}). Hence, without loss of generality, we can assume that
the only positive integer $d$ such that $d|\overline{k}$ and $d|\overline{m}$ is $1$;
which we denote by 
\[ (\overline{k}, \overline{m}) = 1.
\]
Equivalently,
\[ \overline{k}_1\mathbf{Z}+\overline{k}_2\mathbf{Z}+\overline{k}_3\mathbf{Z}+
   \overline{m}_1\mathbf{Z}+\overline{m}_2\mathbf{Z}+\overline{m}_3\mathbf{Z}+
   =\mathbf{Z}.
\]

Finally, suppose there exists some positive integer $d\neq 1$ such that
$d|a$ and $d|b$. Then $d|n$; set
\[ n_0=\frac{1}{d}n,\ \ a_0=\frac{1}{d}a,\ \ b_0=\frac{1}{d}b.
\]
Observe that
\[ G_{i,j}^{ir}(a_0\overline{k}, b_0\overline{m}) = 
    \frac{1}{d^3} G_{i,j}^{ir}(a\overline{k}, b\overline{m}) = 0.
\]
It follows from lemma 3.5 (2) of \cite{B-M-N4} that the vector $n_0$ of the resonant
triplet $(n_0, a_0\overline{k}, b_0\overline{m})$ can also be written
as 
\[ n_0=\hat{k}+\hat{m} \mbox{ with } (n_0, \hat{k}, \hat{m}) \mbox{ verifying (\ref{3wave}).}
\]
But then
\[ n=dn_0 = a\overline{k}+b\overline{m}
    =a^\prime \sigma_i(\overline{k})+b^\prime \sigma_j(\overline{m})
    = d\hat{k}+d\hat{m}.
\]
From lemma 3.5 (2) of \cite{B-M-N4}, $(d\hat{k}, d\hat{m})$ must coincide with either 
one of the pairs 
\[ (a^\prime\sigma_i(\overline{k}), b^\prime \sigma_j(\overline{m})),\ \ 
   (b^\prime \sigma_j(\overline{m}), a^\prime\sigma_i(\overline{k})).
\]
In particular, $d|a^\prime\overline{k}$ and $d|b^\prime\overline{m}$.
Since $d|a$ and $(a, a^\prime )=1$, we have  $(d, a^\prime )$; 
similarly $(d,b^\prime )=1$. But then Euclid's lemma yields that
$d|\overline{k}$ and $d|\overline{m}$, which contradicts the fact that 
$(\overline{k}, \overline{m})=1$. 
Hence we have proven that $(a,b)=1$. In a similar way, one can show that $(a^\prime ,b^\prime )=1$.

{\bf Conclusion:} It follows from the above study that $n\in\mathbf{Z}^*$
admits the two decompositions
\[ n=a\overline{k}+b\overline{m}=
  a^\prime \sigma_i(\overline{k})+b^\prime \sigma_j(\overline{m})
\]
with 
\[ (a, a^\prime ) =  (b, b^\prime ) = (a,b) = (a^\prime , b^\prime ) = 1,\ \ 
   (\overline{k}, \overline{m} ) = 1.
\]

The triplets $(n, a\overline{k}, b\overline{m})$ and 
$(n, a^\prime\sigma_i(\overline{k}), b^\prime \sigma_j(\overline{m}))$ both verify 
the resonant condition (\ref{3wave}) (from the homogeneity of this condition)
as well as the condition of non-catalyticity (\ref{noncat}).
Indeed, $aba^\prime b^\prime \neq 0$ and the condition (\ref{noncat}) on the initial
triplet $(n,k,m)$ imply that the reduced triplet $(n, \overline{k},\overline{m})$ also
verifies (\ref{noncat})). Finally, the degeneracy condition (\ref{dege})
\[ G_{i,j}^{ir}(a\overline{k}, b\overline{m}) = 0
\]
is verified.
%-------------------------------------------------
\medskip

{\bf Acknowledgments.}  We would like to thank A.I. Bobenko,  C. Bardos and G. Seregin
for very useful discussions. The assistance of Dr. B. S. Kim is gratefully
acknowledged. A.M. and B.N. acknowledge  the support of the AFOSR contract FA9550-05-1-0047. 

\medskip

%------------------------------------------------------


\begin{thebibliography}{99}

\bibitem[Ar1]{Ar1}
 Arnold, V.I., Mathematical methods of classical mechanics,
 Springer-Verlag, New York-Berlin, 1978.

\bibitem[Ar2]{Ar2}
Arnold, V.I., Small denominators. I. Mappings of the
circumference onto itself, {\it Amer. Math. Soc. Transl. Ser. 2},
{\bf 46} (1965), p. 213-284.

\bibitem[Ar-Khe]{Ar-Khe}
Arnold, V.I. and  Khesin, B.A.,  Topological Methods in
Hydrodynamics, Applied Mathematical Sciences, {\bf 125}, Springer,  1997.

\bibitem[B-M-N1]{B-M-N1}
Babin, A., Mahalov, A. and  Nicolaenko, B.,  Global splitting, 
integrability and regularity of 3D Euler and Navier-Stokes equations
for uniformly rotating fluids, {\it European J. Mechanics B/Fluids},
{\bf 15}  (1996), p. 291-300.

\bibitem[B-M-N2]{B-M-N2}
Babin, A., Mahalov, A. and  Nicolaenko, B., Global regularity
and integrability of 3D Euler and Navier-Stokes equations for
uniformly rotating fluids, {\it Asymptotic Analysis}, {\bf15}  (1997), p.
103--150.

\bibitem[B-M-N3]{B-M-N3}%{indiana}
Babin, A., Mahalov, A. and  Nicolaenko, B., Global regularity
of 3D rotating Navier-Stokes equations for resonant domains, {\it
Indiana Univ. Math. J.}, {\bf 48}  (1999), No. 3, p. 1133-1176.

\bibitem[B-M-N4]{B-M-N4}%{indiana1}
Babin, A., Mahalov, A. and  Nicolaenko, B., 3D Navier-Stokes
and Euler equations with initial data characterized by uniformly
large vorticity, {\it Indiana Univ. Math. J.}, {\bf 50}  (2001), p. 1-35.


\bibitem[B-K-M]{B-K-M}%{bkm}
Beale, J.T.,  Kato, T. and Majda, A.,  Remarks on the breakdown
of smooth solutions for the 3D Euler equations, {\it Commun. Math.
Phys.}, {\bf 94}  (1984), p. 61-66.


\bibitem[Bes]{Bes}%{besic} 
 Besicovitch, A.S., Almost Periodic Functions,
Dover, New York,  1954.

\bibitem[Bo-Mi]{Bo-Mi}%{BM}
Bogoliubov, N.N. and  Mitropolsky, Y. A.,
Asymptotic Methods in the Theory of Non-linear Oscillations,
Gordon and Breach Science Publishers, New York,  1961.

\bibitem[Bou-Br]{Bou-Br}%{brezis}
Bourguignon, J.P. and Brezis, H., Remark on the Euler
equations, {\it J. Func. Anal.}, {\bf 15}  (1974), p. 341-363.


\bibitem[Ch-Ch-Ey-H]{Ch-Ch-Ey-H}  Chen, Q.,   Chen, S.,   Eyink, G.L. and
 Holm, D.D., Intermittency in the joint cascade of
energy and helicity, {\it Phys. Rev. Letters}, {\bf 90}  (2003),
p. 214503.

\bibitem[Cor]{Cor}%{cord} 
 Corduneanu, C.,
Almost periodic Functions, Wiley-Interscience, New York,  1968.

%\bibitem[DiPe]{DiPe} DiPerna, R.J., Convergence of the viscosity method for
%  isentropic gas dynamics, Comm. Math. Phys., {\bf 91} (1983), p. 1-30.

\bibitem[DiPe-Li]{DiPe-Li}  DiPerna, R.J. and Lions, P.L., Ordinary differential equations,
   Sobolev spaces and transport theory, Invent. Math., {\bf 98} (1989), p. 511-547.

\bibitem[Fe]{Fe} Fefferman, C.L., Existence and smoothness of the Navier-Stokes equations,
  \textit{The millennium prize problems}, Clay Math. Inst., Cambridge, MA (2006), p. 57-67.


\bibitem[Fri]{Fri}Frisch,  U., Turbulence: the legacy of A. N. Kolmogolov,
  Cambridge University Press, 1995.

\bibitem[Fro-M-N]{Fro-M-N} Frolova, E., Mahalov, A. and Nicolaenko, B.,
  Restricted interactions and global regularity of
3D rapidly rotating Navier-Stokes equations in cylindrical domains,  
Journal of Mathematical Sciences,
Springer, to appear. 

\bibitem[Gl1]{Gl1}  Gledzer, E.B., Systema gidrodinamicheskovo
tipa, dopuskayuchaya dva kvadratichnykh integrala dvizheniya,
{\it D. A. N. USSR}, {\bf 209}  (1973), No. 5.


\bibitem[G-D-O]{G-D-O}  Gledzer, E.B.,  Dolzhanski, F.V. and 
   Obukhov, A.M., Systemi gidrodinamitcheskovo tipa i
ikh primetchnii, Nauka, Moscow,  (1987).

\bibitem[G-M-N]{GMNNew}%{bardos}  
Golse, F., Mahalov, A.,  Nicolaenko, B., 
in preparation.

\bibitem[Gu-Ma]{Gu-Ma} Guckenheimer, J. and Mahalov, A., Resonant
triad interaction in symmetric systems, Physica D, {\bf 54} (1992),
267-310.

\bibitem[Hou1]{Hou1} Hou, T.Y., Deng, J. and Yu, X.,
Geometric properties and nonblowup of 3D incompressible Euler flow,
C.P.D.E., {\bf 30} (2005), p. 225-243.

\bibitem[Hou2]{Hou2} Hou, T.Y. and Li, R., Dynamic depletion of vortex
stretching and non-blowup of the 3D incompressible Euler equations,
CALTECH, preprint (2006).


\bibitem[Ka]{Ka}%{kato} 
 Kato, T., Nonstationary flows of viscous and
ideal fluids in $\bf{R}^3$, {\it J. Func. Anal.}, {\bf 9}  (1972), p.
296-305.

\bibitem[Ke]{Ke} Kerr, R.M., Evidence for a singularity of the three
dimensional, incompressible Euler equations, Phys. Fluids, {\bf 5} (1993),
No. 7, p. 1725-1746.


\bibitem[Les]{Les}  Lesieur, M., Turbulence in fluids,
 2nd edition, Kluwer, Dortrecht,  1990.

\bibitem[Li]{Li} Lions, P.L.,
  Mathematical Topics in Fluid Mechanics: Incompressible Models Vol 1,
 Oxford University Press, 1998.

\bibitem[Mah]{Mah} Mahalov, A., The instability of rotating fluid columns
subjected to a weak external Coriolis force, Phys. Fluids A, {\bf 5} (1993),
No. 4,  p. 891-900.


\bibitem[M-N-B-G]{M-N-B-G}%{bardos}  
Mahalov, A.,  Nicolaenko, B., Bardos,  C. and  Golse, F.,
Non blow-up of the 3D Euler equations for a class of three-dimensional
 initial data in cylindrical domains, {\it Methods and Applications of
Analysis}, {\bf 11}  (2004), No. 4, p. 605-634.

\bibitem[Man]{Man} Manakhov, S.V., Note on the integration of
Euler's equations of the dynamics of a $n$-dimensional rigid
body, Funct. Anal. and Appl., {\bf 10} (1976), No. 4, p. 328-329.


\bibitem[Mor1]{Mor1}  Moreau, J.J., Une methode de cinematique
fonctionelle en hydrodynamicque, C.R. Acad. Sci. Paris,
{\bf 249}  (1959), p. 2156-2158

\bibitem[Mor2]{Mor2}  Moreau, J.J., Constantes d'un il\^{o}t
tourbillonaire en fluide parfait barotrope, 
C.R. Acad. Sci. Paris,
{\bf 252}  (1961), p. 2810-2812

\bibitem[Mof]{Mof}  Moffatt, H.K., The degree of knottedness
of tangled vortex lines, {\it J. Fluid Mech.}, {\bf 106}  (1969), p. 117-129.


\bibitem[Poi]{Poi}%{poi}
 Poincar\'{e}, H., Sur la pr\'{e}cession des corps
d\'{e}formables, {\it Bull. Astronomique}, {\bf 27}  (1910), p. 321-356.


\bibitem[Sob]{Sob}%{sobolev}
 Sobolev, S.L., Ob odnoi novoi zadache matematicheskoi
fiziki, {\it Izvestiia Akademii Nauk SSSR, Ser. Matematicheskaia},
{\bf 18}  (1954), No. 1, p. 3--50.

\bibitem[Vish]{Vish}  Vishik, S. M., Ob invariantnyh characteristikah
 kvadratichno-nelineynyh sistem kaskadnovo tipa, {\it D. A. N. USSR}, 
{\bf 228}  (1976), No. 6, p. 1269-1270.

\bibitem[We-Wil]{We-Wil} Weiland, J. and Wilhelmsson, H., Coherent nonlinear
interactions of waves in plasmas, Pergamon, Oxford, 1977.


\bibitem[Yu1]{Yu1}%{yudovich}
 Yudovich, V.I., Non stationary flow of an ideal
incompressible liquid, {\it Zb. Vych. Mat.}, {\bf 3}  (1963), p. 1032-1066

\bibitem[Yu2]{Yu2}%{yudovich}
 Yudovich, V.I., Uniqueness theorem for the basic
 nonstationary problem in th dynamics of an ideal
incompressible fluid, {\it Math. Res. Letters}, {\bf 2}  (1995), p. 27-38.

\bibitem[Zak-Man1]{Zak-Man1} Zakharov, V.E. and Manakov, S.V., Resonant interactions of wave packets in nonlinear
media, Sov. Phys. JETP Lett., {\bf 18} (1973), 243-245.

\bibitem[Zak-Man2]{Zak-Man2} Zakharov, V.E. and Manakov, S.V., The theory of resonance interaction
of wave packets in nonlinear media, Sov. Phys. JETP, {\bf 42} (1976), 842-850.



\end{thebibliography}
\end{document}